\documentclass[preprint,11pt,authoryear]{elsarticle}

\makeatletter
\def\ps@pprintTitle{%
  \let\@oddhead\@empty
  \let\@evenhead\@empty
  \def\@oddfoot{\reset@font\hfil\thepage\hfil}
  \let\@evenfoot\@oddfoot
}
\makeatother

\usepackage[utf8]{inputenc}
\usepackage{geometry}
\usepackage{fullpage,times}
\usepackage{multirow}
\usepackage{comment}
\usepackage{color,soul}
\usepackage{mathrsfs}
\usepackage{amsmath, amsfonts,amssymb}
\usepackage{graphicx}
\usepackage{xspace}
\usepackage{multirow}
\usepackage{floatrow}
\usepackage{eucal}
\usepackage{comment}
\usepackage{mathtools}
\usepackage{pgfgantt}
\usepackage{setspace}
\usepackage{booktabs}
\usepackage[ruled,vlined,linesnumbered,noend]{algorithm2e}
\usepackage{subfig}
\usepackage{enumitem}
\usepackage{booktabs}
\usepackage{makecell}

\usepackage{tikz}
\usetikzlibrary{arrows.meta,calc,fit,positioning}

\usepackage{url}
\usepackage[authoryear]{natbib}

\usepackage{hyperref}
\hypersetup{
    colorlinks=true, %set true if you want colored links
    linkcolor=blue,%blue,  %choose some color if you want links to stand out
    citecolor = blue,
    urlcolor = blue,
    filecolor=blue,
}

\newcommand{\cA}{\ensuremath{\mathcal{A}}}

\newcommand{\cD}{\ensuremath{\mathcal{D}}}

\newcommand{\cG}{\ensuremath{\mathcal{G}}}

\newcommand{\cK}{\ensuremath{\mathcal{K}}}
\newcommand{\cL}{\ensuremath{\mathcal{L}}}

\newcommand{\cP}{\ensuremath{\mathcal{P}}}

\newcommand{\cS}{\ensuremath{\mathcal{S}}}
\newcommand{\cT}{\ensuremath{\mathcal{T}}}

\newcommand{\cV}{\ensuremath{\mathcal{V}}}

\begin{document}
%%%%%%%%%%%%%%%%%%%%%%%%%%%%%%%%%%%%%%%%%%%%%%%%%%%%%%%%%%%%%%%%%%%%%%%%
\begin{frontmatter}

\title{Fairness-aware Strategic Design of Station-based Electric Car-Sharing Systems}

\author[1]{Jue Zhou}
\ead{zhj5226330@gmail.com}

\author[2]{Zoha Sherkat-Masoumi}
\ead{zoha.sherkatmasoumi@mail.utoronto.ca}

\author[1]{Merve Bodur\corref{cor1}}
\ead{merve.bodur@ed.ac.uk}

\cortext[cor1]{Corresponding author}
\address[1]{School of Mathematics and Maxwell Institute for Mathematical Sciences, University of Edinburgh, Edinburgh, UK}
\address[2]{Department of Mechanical and Industrial Engineering, University of Toronto, Toronto, Ontario, Canada}

%%%%%%%%%%%%%%%%%%%%%%%%%%%%%%%%%%%%%%%%%%%%%%%%%%%%%%%%%%%%%%%%%%%%%%%%
%%%%%%%%%%%%%%%%%%%%%%%%%%%%%%%%%%%%%%%%%%%%%%
%%%%%%%%%%%%% Section: Abstract %%%%%%%%%%%%%%
%%%%%%%%%%%%%%%%%%%%%%%%%%%%%%%%%%%%%%%%%%%%%%
\begin{abstract}
Electric car-sharing systems are pivotal for sustainable urban mobility, but their strategic design is complicated by operational constraints, particularly those arising from the charging needs of electric vehicles. The success of these systems hinges on integrating long-term investment decisions (such as station locations, charger capacities, and fleet size) with daily operational realities, including vehicle routing to serve user trip requests and battery management. While existing integrated models address this strategic-operational link, they have prioritized economic efficiency, overlooking the critical dimension of service equity. This paper addresses this gap by making fairness a central design principle, operationalized through two distinct paradigms, namely,  service-rate disparity and max-min fairness, measured explicitly via realized group service rates rather than static spatial accessibility. To capture demand heterogeneity, we adopt a multi-day representative-demand setting, and develop a bi-objective trajectory-based formulation that jointly optimizes revenue and service equity. We develop a solution framework in which a branch-and-price algorithm solves the single-objective variants of the models, embedded within an exact bi-objective procedure to generate the Pareto frontier and complemented by a diving-heuristic-based approach for obtaining high-quality frontier approximations for larger instances.
Through extensive computational experiments, including a Vienna-based real-data case study, we provide key managerial insights into the fundamental trade-offs between revenue, equity, and system design, demonstrating that the proposed framework can serve as a useful decision-support tool for designing station-based electric car-sharing systems that are both economically viable and socially inclusive. 
%\\
%\mb{[Can add some other statistics based on the results]}
\end{abstract}


\end{frontmatter}
%%%%%%%%%%%%%%%%%%%%%%%%%%%%%%%%%%%%%%%%%%%%%%%%%%%%%%%%%%%%%%%%%%%%%%%%
%\newpage
%\section*{Storyline}
%\include{Zoha_Merve_Carsharing_Paper/storyline}
%\newpage

%%%%%%%%%%%%%%%%%%%%%%%%%%%%%%%%%%%%%%%%%%%%%%%%%%%%%%%%%%%%%%%%%%%%%%%%
%%%%%%%%%%%%%%%%%%%%%%%%%%%%%%%%%%%%%%%%%%%%%%%%%%%%%%%%%%%%%%%%%%%%%%%%
\section{Introduction}
\label{sec.introduction}
%%%%%%%%%%%%%%%%%%%%%%%%%%%%%%%%%%%%%%%%%%%%%%%%%%%%%%%%%%%%%%%%%%%%%%%%
%%%%%%%%%%%%%%%%%%%%%%%%%%%%%%%%%%%%%%%%%%%%%%%%%%%%%%%%%%%%%%%%%%%%%%%%

%% HIGHEST-LEVEL MOTIVATION, in particular with an ENVIROMENTAL FOCUS %%
%% BENEFITS of CAR-SHARING, including ELECTRIFICATION %% 
Rapid urbanization and the urgent need to reduce greenhouse gas emissions are transforming how cities plan and operate their transportation systems. With more than half of the global population already living in urban areas and the share projected to reach 68\% by 2050 \citep{owid-urbanization}, providing {\it low-carbon, space-efficient, and socially inclusive mobility solutions} has become a central priority for sustainable urban development. 
Among the strategies supporting this transition, {\it car sharing} has emerged as an effective means of reducing private vehicle ownership and the environmental footprint of urban transport \citep{loose2010state}, while electrification has further amplified these benefits through lower life-cycle emissions \citep{reportICCT}. This trend is reinforced by strong market growth projections, with the global car-sharing market expected to expand from about US\$11 billion in 2025 to nearly US\$29 billion by 2030 and to serve more than 73 million users worldwide, driven in part by the increasing electrification of shared fleets \citep{websiteMordor,websiteStatista,websiteINVERS}. 

Within this broader context, station-based electric vehicle (EV) car-sharing systems have emerged as important components of urban mobility. As these systems move from early pilots to permanent urban services, the focus has expanded beyond pure economic efficiency to include social equity. In this setting, it is no longer sufficient to ask whether an EV car-sharing system is economically viable or operationally efficient; it is also necessary to ask whether it delivers equitable realized service outcomes across different user groups. 
Traditionally, equity in shared mobility has been studied primarily through the lens of \textit{spatial accessibility}, focusing on whether stations and vehicles are geographically distributed so as to cover diverse neighborhoods and demographic groups \citep{meng2021JTG,qian2022equitable}. However, recent studies have increasingly challenged this accessibility-centric view, arguing that the mere physical presence of infrastructure does not guarantee equitable service outcomes \citep{anaya2022}. Even when stations are located within disadvantaged or peripheral communities, users in these areas often face significantly higher rates of unfulfilled requests compared to those in central, high-demand corridors \citep{meng2021JTG,chen2024,jin2024}. This discrepancy reveals a critical gap: spatial accessibility is a static strategic notion that does not capture the dynamic, revenue-driven nature of daily operations.

This gap is particularly important in EV car-sharing systems because the physical constraints of electric propulsion directly shape operational decisions. Unlike internal combustion engine vehicles, EVs are governed by finite battery capacities and time-consuming charging processes \citep{Davatgari2024EJOR}. From an operational standpoint, serving a long-distance trip requires not only a high state-of-charge, but also substantial post-trip downtime for charging, which reduces the vehicle availability for subsequent requests. In a revenue-driven environment, these operational realities create natural incentives to favor short, high-turnover trips over longer or more energy-intensive ones. As a result, a system that appears fair from the standpoint of spatial accessibility may still generate systematically unequal realized service outcomes across different demand groups, creating a systematic disparity in \emph{realized service outcomes} that undermines the promise of inclusive mobility.

These observations suggest that equity should be embedded directly into the strategic planning of EV car-sharing systems, rather than treated only as an ex-post evaluation criterion. Doing so requires an \textit{operations-aware strategic design} perspective that internalizes daily operational dynamics into long-term planning decisions. Yet the literature does not fully provide such a framework: research on equity in shared mobility has mainly emphasized spatial accessibility, while strategic design models for EV car-sharing systems have primarily focused on efficiency. Thus, it remains unclear how realized service fairness can be incorporated into strategic EV car-sharing design under charging and battery constraints.

This paper addresses this gap by establishing \emph{fairness as a central design principle} in the strategic planning of electric car-sharing systems. We focus on station-based, one-way electric car-sharing systems, a configuration that is particularly well suited to coordinated fleet management and the integration of dedicated charging infrastructure \citep{illgen2019literature,huang2020planning}. Our primary conceptual contribution is the systematic formulation and quantitative comparison of two distinct fairness paradigms, namely, service-rate disparity and max-min fairness, within an integrated strategic-operational planning framework. We measure fairness through the \emph{service rate}, defined as the proportion of fulfilled trip requests within each user group, where user groups are defined according to trip energy consumption. To model the resulting fairness-aware planning problem, we consider a \emph{multi-day representative-demand setting}, in which one common strategic design must support multiple heterogeneous operating days, and develop a \emph{bi-objective trajectory-based formulation}. 

The proposed model, while intuitive in its representation of a vehicle’s complete operational path as a trajectory, gives rise to a large-scale bi-objective mixed-integer program. This, in turn, motivates our second major contribution: a corresponding solution framework in which a tailored branch-and-price algorithm solves the single-objective variants of the models, embedded within an exact bi-objective solution approach to generate the Pareto frontier and complemented by a diving-heuristic-based approach for obtaining high-quality frontier approximations for larger instances. 
 
Our computational study demonstrates both the methodological value and the practical insights of the proposed framework across multiple instance classes, spanning synthetic grid-graph instances, real-world Vienna-based instances, and their extensions. Methodologically, the results show that the proposed solution approach is computationally effective across the studied settings. This includes the revenue-maximization case with a single representative day, the variant studied in the literature, for which our branch-and-price approach obtains optimal or near-optimal solutions for substantially larger cases, including instances with significantly more trips and much finer operational-stage time discretization than those tractable by alternative formulations. Practically, our experiments quantify the Pareto efficiency--equity frontier and reveal how different fairness paradigms can fundamentally reshape system design, for example by shifting investments from vehicle procurement toward charging infrastructure to better support high-energy trip groups. These findings provide actionable insights for planners and policymakers seeking shared mobility systems that are operationally feasible, economically viable, and socially inclusive.

The rest of the paper is organized as follows. Section \ref{sec.literature} reviews the relevant literature. Section \ref{sec.description} presents the problem description. Section \ref{sec.formulation} introduces our trajectory-based  formulation. Section \ref{sec.solution} details the proposed solution methodology. Section \ref{sec.experiments} describes the experimental setup and reports the computational results. Finally, Section \ref{sec.conclusion} provides concluding remarks.

%%%%%%%%%%%%%%%%%%%%%%%%%%%%%%%%%%%%%%%%%%%%%%%%%%%%%%%%%%%%%%%%%%%%%%%%
%%%%%%%%%%%%%%%%%%%%%%%%%%%%%%%%%%%%%%%%%%%%%%%%%%%%%%%%%%%%%%%%%%%%%%%%
\section{Literature Review}
\label{sec.literature}
%%%%%%%%%%%%%%%%%%%%%%%%%%%%%%%%%%%%%%%%%%%%%%%%%%%%%%%%%%%%%%%%%%%%%%%%
%%%%%%%%%%%%%%%%%%%%%%%%%%%%%%%%%%%%%%%%%%%%%%%%%%%%%%%%%%%%%%%%%%%%%%%%

To situate our research within the broader context, 
we begin with studies directly related to car-sharing systems (Section \ref{sec.lit_cs}). We then examine the evolution from spatial accessibility to realized service equity in shared mobility to contextualize our primary contribution (Section \ref{sec.lit_fair}). 

\subsection{Strategic design of car-sharing systems}\label{sec.lit_cs}
Strategic planning is a central theme in the car-sharing literature, focusing on long-term, capital-intensive decisions. A large body of works focuses exclusively on station location and fleet sizing \citep{de2012optimization, correia2014added, boyaci2015optimization,li2016design}. 
However, a common limitation of this research is the simplified treatment of day-to-day operational activities. A smaller stream incorporates limited operational considerations—such as relocation needs, temporal demand variation, or uncertainty—but only at an aggregate level, without explicitly modeling the dynamics that determine service performance \citep{huang2018solving, lu2018optimizing, zhang2021optimizing,kosunda2023optimizing}. At the same time, only a small subset of these studies consider demand fluctuations across different scenarios, such as weekdays and weekends \citep{boyaci2015optimization,Clavijo-Lopez2026TRC}.

To address the suboptimal decisions that can arise from the disconnect between planning and implementation, the notion of operations-aware strategic design has been proposed \citep{chen2023}. This paradigm has also proven valuable in micro-mobility systems  \citep{chrisnawati2025}. 
While \citet{Brandstatter2020TS}, \citet{Santos2023TRC}, and \citet{Davatgari2024EJOR} have advanced more realistic strategic planning for electric car-sharing systems via charging station location modeling, space–time–energy flow integration, and explicit energy feasibility-aware capacity design, respectively, these works remain primarily efficiency-oriented and do not account for fairness in realized service outcomes. 

For car-sharing systems, the inclusion of fairness considerations fundamentally reshapes the underlying optimization problem: it requires a shift from the widely adopted single-objective programming framework to multi-objective or bi-objective programming that can balance competing efficiency and fairness goals. Despite this conceptual necessity, this class of fairness-centric multi-objective optimization problems remains under-explored in the existing literature. A small but growing body of work has applied multi-objective formulations to car-sharing system design, but these studies are limited to trade-offs across operational, economic, and environmental dimensions. Representative works include: the bi-objective model from \citet{Enzi2022ORSp} that minimizes total system cost while maximizing user satisfaction; the multi-objective framework from \citet{Miao2019Energy} that maximizes total served trip distance and minimizes total system cost; the bi-objective model from \citet{boyaci2015optimization} that jointly maximizes the operator’s net revenue and the users’ net benefit; and related tri-objective formulations that simultaneously minimize system cost, $CO_2$ emissions, and unmet travel demand.

\subsection{From spatial accessibility to realized service equity}\label{sec.lit_fair}
When strategic and operational decisions are tightly linked, the majority of existing works has not been designed to evaluate realized equity outcomes.
Early equity assessments in shared mobility were grounded primarily in spatial accessibility. The emphasis was placed either on \textit{horizontal equity}, which promotes a more even spatial distribution of infrastructure to ensure comparable geographic access, or on \textit{vertical equity}, which prioritizes resource allocation toward vulnerable populations to compensate for mobility disadvantages \citep{caggiani2020,giuffrida2023JTG,louro2025TRD}.
As social exclusion has received increasing attention in transportation systems, researchers have also begun to examine the broader barriers faced by disadvantaged groups, such as low-income and elderly populations, in accessing and using shared-mobility services. This broader perspective has contributed to growing interest in mobility justice as a framework for evaluating equity \citep{henriksson2022,haxhija2024RTBM}.

More recent studies, however, make clear that physical coverage alone does not guarantee equitable service outcomes \citep{anaya2022,hosford2024TRD}.
Empirical studies indicate that evenwhen infrastructure is located in disadvantaged neighborhoods, substantial disparities in realized service outcomes can persist because of economic barriers, digital divides, or imbalanced system dispatching \citep{meng2021JTG,chen2021JTG, jin2024}. In response, equity metrics have increasingly moved beyond static spatial measures (e.g., Gini coefficients, spatial coverage rates) \citep{duran2021,cheng2024} toward dynamic, service-quality-oriented indicators, including time-varying service levels \citep{duran2021}, service fulfillment risks \citep{chen2023}, and realized service rates \citep{wang2021fairness}. 

Notably, most of this work on dynamic equity metrics has been concentrated in micromobility systems \citep{duran2021, chen2023, yang2025}. By contrast, similar research tailored to car-sharing remains limited, with only a few studies beginning to examine disparities in service acquisition probabilities and time costs induced by system operational mechanisms \citep{wang2021fairness}. This suggests that, although the literature has started to move from spatial accessibility toward realized service outcomes, that transition remains underdeveloped in the car-sharing context.

\subsection{Research gap}
Although the literature has made substantial progress on transportation equity, electrification-aware modeling, and integrated strategic-operational planning, an important gap remains at the intersection of these streams. Research on equity in shared mobility has increasingly moved toward dynamic, outcome-based notions of fairness, but this development has been concentrated largely in micromobility systems. At the same time, the literature on electric car-sharing and operations-aware strategic design has made substantial progress in modeling charging constraints, fleet operations, and infrastructure planning, but has remained focused primarily on profit maximization or cost minimization.

To the best of our knowledge, the literature still lacks an operations-aware strategic planning framework for station-based one-way EV car-sharing that embeds realized-service fairness directly into long-term design while preserving battery and charging feasibility across multiple representative days. This gap motivates the present study, which develops a fairness-aware, operations-aware framework for jointly modeling infrastructure decisions, vehicle charging dynamics, and realized service equity.

%%%%%%%%%%%%%%%%%%%%%%%%%%%%%%%%%%%%%%%%%%%%%%%%%%%%%%%%%%%%%%%%%%%%%%%%
%%%%%%%%%%%%%%%%%%%%%%%%%%%%%%%%%%%%%%%%%%%%%%%%%%%%%%%%%%%%%%%%%%%%%%%%
\section{Problem Description}
\label{sec.description}
%%%%%%%%%%%%%%%%%%%%%%%%%%%%%%%%%%%%%%%%%%%%%%%%%%%%%%%%%%%%%%%%%%%%%%%%
%%%%%%%%%%%%%%%%%%%%%%%%%%%%%%%%%%%%%%%%%%%%%%%%%%%%%%%%%%%%%%%%%%%%%%%%

We study the strategic design of a one-way, station-based electric car-sharing system under explicit fairness considerations. The planner selects a long-term infrastructure and fleet configuration and evaluates that design through its induced operational performance over a set of representative operating days. The planning problem therefore combines three elements: infrastructure and fleet design, day-level operational feasibility, and fairness in realized service outcomes across user groups. The main notation used throughout this paper is summarized in Table~\ref{tab:problem_inputs}.

\begin{table}[h!]
\small
\centering
\begin{tabular}{p{1.5cm}p{12.5cm}}
\toprule
Notation & Description \\
\midrule
$\cS$ & Set of candidate station locations \\
$\cD$ & Set of representative operating days \\
$\cG$ & Set of user groups used to evaluate service equity \\
$\pi_d$ & Normalized weight (occurrence probability) associated with representative day $d \in \cD$ \\
$\cT_d$ & Set of discrete time periods for representative day $d \in \cD$ \\
$\cK_d$ & Set of trip requests observed on representative day $d \in \cD$ \\
$\cK^g_{d}$ & Set of trip requests associated with user group $g \in \cG$ on representative day $d \in \cD$ \\
$s_k,e_k$ & Start time and end time of trip request $k$ \\
$\cS_k^+,\cS_k^-$ & Accessible pick-up station set and accessible drop-off station set for trip $k$ \\
$b_k$ & Driving-energy requirement of trip $k$ \\
$r_k$ & Revenue collected if trip $k$ is served \\
$\cP_d$ & Set of feasible vehicle trajectories for representative day $d \in \cD$ \\
\bottomrule
\end{tabular}
\caption{Input data and system definition}
\label{tab:problem_inputs}
\end{table}

Let $\cS$ denote the set of candidate station locations and $\cD$ the set of representative operating days. Each representative day $d \in \cD$ is associated with a weight $\pi_d$, a time-discretized planning horizon $\cT_d$, and a realized set of trip requests $\cK_d$. Strategic decisions are taken once for the entire system, while their consequences are evaluated over all representative days. The quality of a candidate design is therefore determined by its ability to support feasible operations across heterogeneous daily demand patterns. We assume that the representative-day weights are normalized so that $\sum_{d\in\cD}\pi_d=1$.

Fairness in this paper is defined as \emph{service access equity} across user groups. The relevant outcome is the service rate achieved by each group, namely the proportion of trip requests from that group that can actually be served. Accordingly, let $\cG$ denote the set of user groups and let $\cK^g_{d} \subseteq \cK_d$ denote the subset of trip requests associated with group $g \in \cG$ on day $d \in \cD$. In the empirical implementation, group membership is inferred from the energy requirement of each requested trip, so that the analysis captures whether more energy-intensive travel needs are systematically disadvantaged. The modeling structure itself is not restricted to two groups and remains valid for a general partition of users.

For each user group, the service rate on a representative day is the proportion of that group's requests that are successfully served on that day. The corresponding average group service rate is then obtained by aggregating these day-level service rates across representative days using the weights $\pi_d$. The fairness paradigms studied in this paper are based on these group-level service outcomes: one protects the worst-served group, while the other controls the disparity in service rates across groups. 
For convenience, we consider $|\cK_d^g|>0$ for every $g\in\mathcal G$ and $d\in\mathcal D$, so that these group-level service rates are well-defined.

The strategic design and daily operations are linked through feasible \emph{vehicle trajectories}. For each representative day $d \in \cD$, let $\cP_d$ denote the set of feasible trajectories. A trajectory represents the complete daily schedule of a single electric vehicle and consists of an ordered sequence of trip events and charging events. Each trip request $k \in \cK_d$ is characterized by a start time $s_k$, an end time $e_k$, a set of accessible pick-up stations $\cS_k^+$, a set of accessible drop-off stations $\cS_k^-$, a driving-energy requirement $b_k$, and an associated revenue $r_k$. A trip event in a trajectory specifies which request is served, which feasible pick-up and drop-off stations are used, when the service occurs, and the associated battery consumption.

\begin{figure}[h]
 \centering
\includegraphics[scale=0.65]{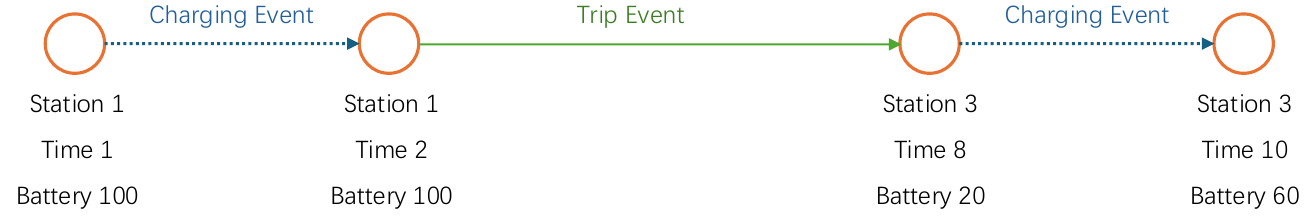}
 \caption{An example of vehicle trajectory}
 \label{sample_trajectory}
\end{figure}

A charging event occurs while the vehicle is parked at a station, either at the beginning of the day or after completing a previous trip, and lasts until the next departure or the end of the planning horizon. Figure~\ref{sample_trajectory} illustrates this trajectory concept for one representative day. The vehicle begins the day at Station~1, departs at time~2 to serve a trip to Station~3, arrives at time~8 with a reduced battery level, and then charges again until the end of the planning horizon.

Before turning to the optimization model, we state the main assumptions underlying this trajectory representation.
\begin{itemize}
    \item The operator procures a homogeneous fleet of EVs. This is standard in strategic design models and avoids confounding infrastructure decisions with vehicle heterogeneity \citep{xu2018electric,xu2021electric}.
    \item All vehicles are fully charged and parked at their assigned stations at the beginning of each representative day. This is a planning-level simplification adopted for tractability \citep{li2016design,boyaci2015optimization,zhao2018integrated,xu2019fleet}. Related work has considered more flexible battery-state modeling, including route-specific initial conditions and partial recharging in fixed-route settings \citep{Davatgari2024EJOR} and, in other settings, inter-period state-of-charge continuity \citep{hellem2021dynamic}.
    %Existing work has relaxed such assumptions to accommodate flexible initial SOC, inter-day SOC continuity, and generalized partial charging \citep{hellem2021dynamic,Davatgari2024EJOR}.
    \item Vehicles are assumed to recharge at a constant, linear rate while parked, which is a common planning-level approximation in related electrified-transportation models \citep{Brandstatter2020TS,huang2020planning,Davatgari2024EJOR}. On the other hand, no such assumption is made for battery discharge during service; instead, discharge is represented by the trip-specific energy requirement recorded in the input data.
    %Battery discharge during a served trip is represented by the trip-specific energy requirement recorded in the input data, and vehicles recharge at a constant, linear rate while parked. This is a common planning-level approximation, consistent with assumptions used in related electrified-transportation models \citep{brandstatter2017determining,huang2020planning,Davatgari2024EJOR}.
    \item Trips are treated as atomic service tasks: once a vehicle starts serving a user, the trip must be completed without interruption, and no en-route charging is allowed while the user is on board. Accordingly, any service option retained in the input data must be executable within a single battery discharge cycle \citep{boyaci2015optimization}.
\end{itemize}

Under these assumptions, each trajectory is individually feasible in terms of timing, location consistency, and battery evolution. Charging events embedded in a trajectory indicate where and when charging capacity would be required if that trajectory were selected, but they do not assume that charging infrastructure has already been installed at those locations. System-wide compatibility with the strategic design is enforced later in the optimization model through station-opening and charger-capacity constraints. This distinction separates vehicle-level schedule feasibility from network-level infrastructure feasibility.

%%%%%%%%%%%%%%%%%%%%%%%%%%%%%%%%%%%%%%%%%%%%%%%%%%%%%%%%%%%%%%%%%%%%%%%%
%%%%%%%%%%%%%%%%%%%%%%%%%%%%%%%%%%%%%%%%%%%%%%%%%%%%%%%%%%%%%%%%%%%%%%%%
\section{Mathematical Formulation}
\label{sec.formulation}
%%%%%%%%%%%%%%%%%%%%%%%%%%%%%%%%%%%%%%%%%%%%%%%%%%%%%%%%%%%%%%%%%%%%%%%%
%%%%%%%%%%%%%%%%%%%%%%%%%%%%%%%%%%%%%%%%%%%%%%%%%%%%%%%%%%%%%%%%%%%%%%%%
This section formulates the strategic planning problem as a trajectory-based mixed-integer program over representative operating days. We first introduce the bi-objective planning model with revenue and fairness objectives. We then specify two alternative fairness functions. Finally, we present the practically relevant single-objective variants used for computation and policy analysis.

\subsection{Bi-objective model}
\label{sec.basic_formulation}
For each representative day $d \in \cD$, the set $\cP_d$ contains all trajectories that are individually feasible with respect to time, location consistency, and battery evolution. Importantly, feasibility at the trajectory level does not yet imply that the corresponding charging infrastructure has been installed. Instead, the charging events embedded in a trajectory indicate the charger occupancy that would be required if that trajectory were selected. The optimization model then decides which day-specific trajectory sets can be jointly supported by the strategic design.
Table~\ref{table:variables} summarizes the model-specific parameters and decision variables. 

\begin{table}[h!]
\small
\centering
\begin{tabular}{p{1cm}p{13.4cm}}
\toprule
\multicolumn{2}{l}{\underline{\sc Model parameters}} \\
$c_s$ & Cost of opening station $s \in \cS$ \\
$c'_s$ & Cost of installing one charger at station $s \in \cS$ \\
$Q_s$ & Maximum charger capacity at station $s \in \cS$ \\
$H$ & Upper bound on the number of vehicles that can be procured \\
$c_\texttt{h}$ & Cost of procuring one vehicle \\
$W$ & Total available investment budget \\
$\hat{a}^{k}_{pd}$ & Equals $1$ if trip $k \in \cK_d$ is served in trajectory $p \in \cP_d$, and $0$ otherwise \\
$\hat{b}^{st}_{pd}$ & Equals $1$ if trajectory $p \in \cP_d$ occupies a charger at station $s \in \cS$ at time $t \in \cT_d$, and $0$ otherwise \\
%$M$ & Sufficiently large penalty coefficient \\
\midrule
\multicolumn{2}{l}{\underline{\sc Decision variables}} \\
$x^d_k$ & Equals $1$ if trip $k \in \cK_d$ is accepted on representative day $d \in \cD$, and $0$ otherwise \\
$y_s$ & Equals $1$ if station $s \in \cS$ is opened, and $0$ otherwise \\
$z_s$ & Number of chargers installed at station $s \in \cS$ \\
$v_\texttt{h}$ & Number of vehicles procured for the fleet \\
$\lambda^d_{p}$ & Equals $1$ if trajectory $p \in \cP_d$ is selected on representative day $d \in \cD$, and $0$ otherwise \\
$u_g^d$ & Service rate of user group $g \in \cG$ on representative day $d \in \cD$ \\
$\bar{u}_g$ & Weighted average service rate of user group $g \in \cG$ across representative days \\
$\eta$ & Lower bound on the average service rate across user groups in the max-min model \\
$\delta$ & Upper bound on the disparity of average service rates in the disparity model \\
$\bar{\eta}$ & Slack variable for the minimum-service guarantee \\
$\bar{\delta}$ & Slack variable for the maximum-disparity requirement \\
\bottomrule
\end{tabular}
\caption{Model parameters and decision variables}
\label{table:variables}
\end{table}

Let $R(x)$ denote the expected revenue over representative days, i.e., 
\[
R(x) = \sum_{d\in\cD}\pi_d\sum_{k\in\cK_d} r_k x_k^d.
\]
Let $F(x)$ denote a generic fairness objective defined on realized group service outcomes; we specify the two fairness definitions considered in Section~\ref{sec.fair_formulation}. 
The bi-objective planning model \textbf{BM} is then written as
\begingroup
\allowdisplaybreaks
\begin{subequations}\label{eq.biobjective}
\begin{align}
\textbf{(BM):} \ \max \  & \left\{R(x),F(x)\right\} \label{eq.biobjective_obj}\\
\hbox{s.t.} \ 
&\sum_{p \in \cP_d} \hat{a}^{k}_{pd} \lambda^d_p = x^d_k && \forall k \in \cK_d,\ \forall d \in \cD \label{eq.trip_coverage}\\
& \sum_{p \in \cP_d} \hat{b}^{st}_{pd} \lambda^d_p \le z_s && \forall s \in \cS,\ \forall t \in \cT_d,\ \forall d \in \cD \label{eq.charger_available}\\
& z_s \le Q_s y_s && \forall s \in \cS \label{eq.charger_capacity}\\
& \sum_{s \in \cS}(c_s y_s + c'_s z_s) + c_\texttt{h} v_\texttt{h} \le W \label{eq.budget}\\
& \sum_{p \in \cP_d} \lambda^d_p \le v_\texttt{h} && \forall d \in \cD \label{eq.vehicle_num}\\
& x^d_k \in \{0,1\} && \forall k \in \cK_d,\ \forall d \in \cD \label{eq.bound_x}\\
& y_s \in \{0,1\} && \forall s \in \cS \label{eq.bound_y}\\
& z_s \in \{0,1,\dots,Q_s\} && \forall s \in \cS \label{eq.bound_z}\\
& v_\texttt{h} \in \{0,1,\dots,H\} \label{eq.bound_h}\\
& \lambda^d_p \in \{0,1\} && \forall p \in \cP_d,\ \forall d \in \cD \label{eq.bound_lambda}
\end{align}
\end{subequations}
\endgroup
Objective~(\ref{eq.biobjective_obj}) aims to maximize expected revenue and fairness simultaneously. 
Constraints~(\ref{eq.trip_coverage}) link trip-acceptance decisions and trajectory-selection decisions on each representative day. Constraints~(\ref{eq.charger_available}) enforce charger availability over time and provide the main strategic-operational linking mechanism because they guarantee that the charging requirements implied by selected trajectories can be supported by the installed infrastructure. Constraints~(\ref{eq.charger_capacity}) restrict charger installation to opened stations and respect the site-specific capacity limits. Constraint~(\ref{eq.budget}) imposes the total investment budget, while constraint~(\ref{eq.vehicle_num}) limits the number of active trajectories on each representative day by the procured fleet size. Finally, constraints~(\ref{eq.bound_x})--(\ref{eq.bound_lambda}) specify the variable domains.

We note that a more explicit formulation could in principle be written using vehicle-indexed operational variables for trip assignment, charging, and battery-state evolution over time. We adopt a trajectory-based formulation because, in the present setting, a trajectory is the natural object representing a complete daily vehicle schedule and therefore provides a direct and compact way to model the strategic-operational coupling. This representation captures timing, location consistency, battery feasibility, and charging requirements within a single object, while making clear that charging events encode required charger occupancy rather than preinstalled infrastructure. The same structure is later exploited by our proposed solution approach.

\subsection{Definitions of fairness objectives}
\label{sec.fair_formulation}
To instantiate the fairness objective, $F(x)$, let $u^d_g$ denote the service rate of user group $g \in \cG$ on representative day $d \in \cD$, and let $\bar{u}_g$ denote the corresponding weighted average service rate across representative days. These quantities are linked to the trip-acceptance variables by
\begingroup
\allowdisplaybreaks
\begin{subequations}\label{eq.service_rate}
\begin{align}
u^d_g &= \frac{1}{|\cK^d_{g}|}\sum_{k \in \cK_{g}^d} x^d_k && \forall g \in \cG,\ \forall d \in \cD,\label{eq.problem_service_rate_day}\\
\bar{u}_g &= \sum_{d \in \cD}\pi_d u^d_g && \forall g \in \cG.
\label{eq.problem_service_rate_avg}
\end{align}
\end{subequations}
\endgroup
The two fairness schemes studied in this paper correspond to two different specifications of $F(x)$.

\subsubsection{Max-min fairness}
\label{sec.fair_formulation_maxmin}
The first fairness concept aims to protect the most disadvantaged user group. We therefore define
\begin{equation}
F^{\mathrm{MM}}(x)=\min_{g\in\cG}\bar{u}_g.
\label{eq.fairness_maxmin}
\end{equation}
This specification follows the classical max-min fairness principle \citep{kumar2006fairness} and focuses attention on the worst-served group.

We introduce a continuous variable $\eta \in [0,1]$ representing the minimum average group service rate and add the linking constraints
\begin{equation}
\eta \le \bar{u}_g, \quad \forall g \in \cG.
\label{eq.maxmin_eta_def}
\end{equation}
Under this representation, maximizing $F^{\mathrm{MM}}(x)$ is equivalent to maximizing $\eta$ subject to~(\ref{eq.maxmin_eta_def}).

\subsubsection{Service-rate disparity}
\label{sec.fair_formulation_parity}
The second fairness concept directly limits the disparity of realized service rates across groups. We define
\begin{equation}
F^{\mathrm{SD}}(x)=-\max_{\substack{g,g'\in\cG\\g\neq g'}} |\bar{u}_g-\bar{u}_{g'}|,
\label{eq.fairness_disparity}
\end{equation}
so that larger values of $F^{\mathrm{SD}}(x)$ correspond to smaller inter-group gaps.

Let $\delta \in [0,1]$ denote the maximum allowable difference in average service rates. For any two groups $g,g' \in \cG$ with $g \neq g'$, we require
\begin{equation}
|\bar{u}_g - \bar{u}_{g'}| \le \delta.
\label{eq.parity_abs}
\end{equation}
For the practically relevant case of two user groups, say a lower-energy group $g_1$ and a higher-energy group $g_2$, the absolute-value constraint is linearized as
\begingroup
\allowdisplaybreaks
\begin{subequations}\label{eq.service_disparity}
\begin{align}
\bar{u}_{g_1} - \bar{u}_{g_2} &\le \delta \label{eq.parity_linear_upper}\\
\bar{u}_{g_1} - \bar{u}_{g_2} &\ge -\delta. \label{eq.parity_linear_lower}
\end{align}
\end{subequations}
\endgroup
Under this representation, maximizing $F^{\mathrm{SD}}(x)$ is equivalent to minimizing $\delta$ subject to~(\ref{eq.parity_abs}) or, in the two-group case studied computationally, subject to~(\ref{eq.parity_linear_upper})--(\ref{eq.parity_linear_lower}).

\subsection{Practically relevant single-objective variants}
\label{sec.single_objective_variants}
The bi-objective formulation~(\ref{eq.biobjective}) defines the overall planning problem, whereas computation and managerial interpretation rely on objective-constrained single-objective variants. In this setting, stringent fairness targets may render the model infeasible under given fairness metrics. To avoid excluding such cases, we introduce slack variables in the fairness constraints and penalize them in the objective with sufficiently large coefficients. This yields relaxed formulations that remain solvable even when the prescribed fairness requirement cannot be fully attained, while ensuring that any violation is used only as a last resort.

For max-min fairness, the planner specifies a minimum acceptable group service guarantee $\epsilon^{\min}$ and solves
\begin{align*}
\textbf{(M2): \ \ }\max \quad & R(x)-M\bar{\eta}\\
\hbox{s.t.\ \ }
& (\ref{eq.problem_service_rate_day}),(\ref{eq.problem_service_rate_avg}),(\ref{eq.trip_coverage})-(\ref{eq.bound_lambda}), (\ref{eq.maxmin_eta_def}) \\
& \eta + \bar{\eta}\ge \epsilon^{\min} \\
& \bar{\eta}\ge 0.
\end{align*}
The slack variable $\bar{\eta}$ provides a relaxed version of the fairness requirement and protects the model from infeasibility when the requested guarantee is too ambitious. A sufficiently large penalty $M$ ensures that the relaxation is used only when necessary.

For service-rate disparity control, the planner specifies a maximum admissible gap $\epsilon^{\text{gap}}$ and solves
\begin{align*}
\textbf{(M3): \ \ }\max \quad & R(x)-M\bar{\delta}\\
\hbox{s.t.\ \ }
& (\ref{eq.problem_service_rate_day}),(\ref{eq.problem_service_rate_avg}),(\ref{eq.trip_coverage})-(\ref{eq.bound_lambda}), (\ref{eq.parity_linear_upper})-(\ref{eq.parity_linear_lower}) \\
& \delta \le \epsilon^{\text{gap}}+\bar{\delta}\\
& \bar{\delta}\ge 0.
\end{align*}
Here, the slack variable $\bar{\delta}$ yields the corresponding relaxed variant and again serves as a hedge against infeasibility. 

By varying $\epsilon^{\min}$ or $\epsilon^{\text{gap}}$, the planner can trace the practically relevant segments of the Pareto frontier while keeping the interpretation of the fairness target transparent. If all fairness-related constraints are removed from these objective-constrained formulations, the model reduces to the pure revenue-maximization formulation, which we denote by  \textbf{(M1)}. This model coincides with  the benchmark revenue-focused formulation in the literature \citep{Brandstatter2020TS} and serves as the revenue anchor of the Pareto frontier. The fairness-oriented anchor solutions used later in the bi-objective procedure correspond to maximizing $\eta$ in the max-min model and minimizing $\delta$ in the disparity model.

%%%%%%%%%%%%%%%%%%%%%%%%%%%%%%%%%%%%%%%%%%%%%%%%%%%%%%%%%%%%%%%%%%%%%%%%
%%%%%%%%%%%%%%%%%%%%%%%%%%%%%%%%%%%%%%%%%%%%%%%%%%%%%%%%%%%%%%%%%%%%%%%%
\section{Solution Methodology}
\label{sec.solution}
%%%%%%%%%%%%%%%%%%%%%%%%%%%%%%%%%%%%%%%%%%%%%%%%%%%%%%%%%%%%%%%%%%%%%%%%
%%%%%%%%%%%%%%%%%%%%%%%%%%%%%%%%%%%%%%%%%%%%%%%%%%%%%%%%%%%%%%%%%%%%%%%%

This section describes the criterion-space framework used to generate the Pareto frontier, the exact branch-and-price algorithm used to solve the underlying oracle problems, and its diving-based heuristic variant.

\subsection{General framework to generate the Pareto frontier}
\label{sec.pareto_exact}

We generate the Pareto frontier in criterion space by the balanced box method of \citet{boland2015criterion}. The method is applied separately to the two fairness paradigms studied in this paper. For a fixed fairness paradigm, let $\Omega$ denote the feasible set of the corresponding planning model, and let
\[
(F(\vartheta),R(\vartheta))
\]
denote the fairness--revenue criterion pair associated with a feasible solution $\vartheta \in \Omega$, where both objectives are to be maximized.

The method maintains a set of already identified nondominated criterion points together with a set of active boxes in the $(F,R)$ plane. Each active box is defined by two known nondominated points
\[
q^\texttt{L}=(\phi^\texttt{L},\rho^\texttt{L}), \qquad q^\texttt{R}=(\phi^\texttt{R},\rho^\texttt{R}),
\]
with $\phi^\texttt{L}<\phi^\texttt{R}$ and $\rho^\texttt{L}>\rho^\texttt{R}$. Thus, any undiscovered nondominated point between these two known points must lie in the rectangle
\[
\mathbb{B}(q^\texttt{L},q^\texttt{R})
:=
\{(\phi,\rho): \phi^\texttt{L}\le \phi \le \phi^\texttt{R},\ \rho^\texttt{R}\le \rho \le \rho^\texttt{L}\}.
\]

The initial box is determined by two anchor points, obtained by solving lexicographical maximization problems: 
\begin{alignat*}{2}
\text{The revenue-first anchor: } & q^{\texttt{rev}}=(F(\vartheta^{\texttt{rev}}),R(\vartheta^{\texttt{rev}})),
\ \ &&
\vartheta^{\texttt{rev}} \in \text{lex}\max\{(R(\vartheta),F(\vartheta)):\vartheta\in\Omega\} \\
\text{The fairness-first anchor: } & q^{\texttt{fair}}=(F(\vartheta^{\texttt{fair}}),R(\vartheta^{\texttt{fair}})),
\ \ &&
\vartheta^{\texttt{fair}} \in \text{lex}\max\{(F(\vartheta),R(\vartheta)):\vartheta\in\Omega\}.
\end{alignat*}
These two points define the initial search box. Throughout the algorithm, active boxes are processed in nonincreasing order of area so that good frontier approximations are obtained early.

Consider an active box $\mathbb{B}(q^\texttt{L},q^\texttt{R})$, and let $m:=(\rho^\texttt{L}+\rho^\texttt{R})/2$
be the midpoint of the box in the revenue dimension. The box is explored in two steps. (1) We search the upper half of the box by solving
\[
\bar q^\texttt{U}\in
\text{lex}\max
\Bigl\{
(F(\vartheta),R(\vartheta)):
\vartheta\in\Omega,\ 
\phi^\texttt{L}\le F(\vartheta)\le \phi^\texttt{R},\
m\le R(\vartheta)\le \rho^\texttt{L}
\Bigr\}.
\]
If this lexicographic search returns the corner point $q^\texttt{L}$, then the upper half contains no additional nondominated point. Otherwise, it returns a new nondominated point $
\bar q^\texttt{U}=(\bar\phi^\texttt{U},\bar\rho^\texttt{U})$. 
(2) The point $\bar q^\texttt{U}$ dominates every point in the lower half of the box whose fairness value is at most $\bar\phi^\texttt{U}$. Therefore, only the remaining lower region needs to be searched. We do so by solving
\[
\bar q^\texttt{L}\in
\text{lex}\max
\Bigl\{
(R(\vartheta),F(\vartheta)):
\vartheta\in\Omega,\ 
\bar\phi^\texttt{U}<F(\vartheta)\le \phi^\texttt{R},\
\rho^\texttt{R}\le R(\vartheta)\le m
\Bigr\}.
\]
If this search returns the corner point $q^\texttt{R}$, then the remaining lower region contains no additional nondominated point. Otherwise, it returns a new nondominated point $
\bar q^\texttt{L}=(\bar\phi^\texttt{L},\bar\rho^\texttt{L})$. 
By construction, any undiscovered nondominated point from the original box must then lie in one of the two smaller boxes
\[
\mathbb{B}(q^\texttt{L},\bar q^\texttt{U})
\qquad\text{and}\qquad
\mathbb{B}(\bar q^\texttt{L},q^\texttt{R}),
\]
which are inserted into the active set. Repeating this process until no active box remains yields the complete Pareto frontier.

Each lexicographic search above is implemented by solving two box-constrained single-objective oracle problems in sequence. For example, $\text{lex}\max(F,R)$ first maximizes fairness over the current box and then maximizes revenue over the same box while fixing fairness at its optimal value; $\text{lex}\max(R,F)$ is defined analogously. In our setting, the revenue-first oracle problems coincide with the objective-constrained models introduced in Section~\ref{sec.single_objective_variants} once the box-induced bounds are added, while the fairness-first oracle problems are obtained analogously by reversing the roles of revenue and fairness.

Lastly, we note that although the balanced box method is developed for biobjective integer programs, its use remains valid in our setting. The method does not rely on convexity of the feasible region; it requires only that the box-induced oracle problems be solved exactly and that the attainable criterion set be finite. In our models, both revenue and fairness values are ultimately induced by finitely many discrete strategic, trip-acceptance, and trajectory-selection decisions, while the continuous fairness variables are auxiliary quantities determined by those decisions and the fixed input data. Therefore, the attainable criterion set is finite, and exact oracle solutions imply that the balanced box method recovers the complete Pareto frontier for the corresponding fairness paradigm.

\subsection{B\&P to solve the single-objective variants}
\label{sec.bap_exact}

The box-constrained single-objective oracle problems arising in Section~\ref{sec.pareto_exact} inherit the trajectory-based structure of the planning formulation. Since the number of feasible trajectories is exponential, direct enumeration of all columns is intractable. We therefore solve these oracle problems by branch-and-price (B\&P), which provides the exact optimization engine required by the balanced box method. The same B\&P framework can also be used directly to solve the practically relevant single-objective formulations introduced in Section~\ref{sec.single_objective_variants}; within the balanced-box method, the oracle problems are box-constrained versions of this same trajectory-based structure.

Only the day-specific trajectory variables $\lambda_p^d$ are generated dynamically. The remaining variables, namely the shared strategic-design variables $(y_s,z_s)$ and $v_\texttt{h}$, the day-specific trip-acceptance variables $x_k^d$, and, when fairness is active, the auxiliary fairness variables, remain explicit in the master problem. Because trajectories are indexed by representative day, column generation decomposes naturally into one pricing subproblem for each day $d \in \cD$ at every iteration.

\subsubsection{Restricted master problem and pricing subproblem}
For any box-induced oracle problem from Section~\ref{sec.pareto_exact}, column generation starts from a restricted master problem (RMP) defined on current subsets of trajectories. For each representative day $d \in \cD$, let $\cP_d' \subseteq \cP_d$ denote the current set of generated trajectories. Restricting the trajectory variables to $p \in \cP_d'$ and relaxing integrality yields the LP relaxation of the current oracle problem. The box-induced objective and bound structure depend on the lexicographic search being performed, whereas the trajectory-based linking constraints remain unchanged.

Let $\alpha_k^d$ (free), $\beta_{st}^d \ (\geq 0)$, and $\mu_d \ (\geq 0)$ denote the dual values associated with the trip-coverage constraints~\eqref{eq.trip_coverage}, charger-availability constraints~\eqref{eq.charger_available}, and day-specific vehicle constraints~\eqref{eq.vehicle_num}, respectively. Then, for a trajectory $p \in \cP_d$ generated on representative day $d$, its reduced cost is
\begin{equation}
rc_p^d = -
\left(
\sum_{k \in \cK_d} \hat a_{pd}^{k}\alpha_k^d
+
\sum_{s \in \cS}\sum_{t \in \cT_d}\hat b_{pd}^{st}\beta_{st}^d
+
\mu_d
\right).
%rc_p = -\left(\sum_{k \in \cK_d} \hat{a}^k_p \alpha_k + \sum_{s \in \cS} \sum_{t\in (\cT_s^-\cup \cT_s^0)} \hat{b}^{st}_p \beta_{st} + \mu_d \right).
\label{eq:reduced_cost}
\end{equation}
When fairness is activated, the additional fairness terms remain in the master problem through the explicit trip-acceptance variables $x_k^d$ and service rate plus fairness measure related variables (e.g., $u_g^d$,  $\bar{u}_g$, $\eta$, and $\delta$),  and their associated constraints. The pricing structure is unaffected: for each representative day, the pricing subproblem still searches for a feasible trajectory with positive reduced cost. Also, importantly, the pricing subproblem can be formulated as a resource-constrained longest path problem. The same pricing structure is used for both revenue-first and fairness-first lexicographic searches; only the master-problem objective and the box-induced bounds change.

At each column-generation iteration, we therefore solve one pricing subproblem per $d \in \cD$. Any trajectory with positive reduced cost is added to the corresponding day-specific column pool $\cP_d'$. Column generation terminates when no representative day yields a positive reduced-cost trajectory, at which point the current RMP solution is optimal for the LP relaxation of the selected single-objective model.

To generate trajectories for day $d$, we use a day-specific time-expanded network as in \citep{Brandstatter2020TS,wu2022multi}. 
Let $\cG_d=(\cV_d,\cA_d)$ denote the network for representative day $d$, where $\cV_d$ and $\cA_d$ denote the vertex set and arc set, respectively. Each vertex represents a station--time pair $v=(s,\tau)$. For each station $s \in \cS$, let $\cT_{sd}^+$ and $\cT_{sd}^-$ denote the sets of pivotal times on day $d$ at which trips depart from and arrive at station $s$, respectively, and let $\cT_{sd}=\cT_{sd}^+\cup\cT_{sd}^-$. We also add a dummy source node $\bar{o}$ and a dummy sink node $\bar{d}$. Thus,
$
\cV_d=\{\bar{o},\bar{d}\}\cup\{(s,\tau)\mid s\in\cS,\ \tau\in\cT_{sd}\}$.

The arc set $\cA_d$ contains four types of arcs. Initialization arcs ($\cA_{\bar{o}}$) connect the source node to the earliest pivotal-time vertex at each station and represent the start of a daily vehicle schedule. Termination arcs ($\cA_{\bar{d}}$) connect the latest pivotal-time vertex at each station to the sink node. Waiting arcs ($\cA_{\texttt{wait}}$) connect consecutive pivotal-time vertices at the same station and represent periods during which a vehicle remains parked and can recharge. Trip arcs ($\cA_{\texttt{trip}}$) represent the service of trip requests: for each trip $k \in \cK_d$ and each feasible pair $(s^+,s^-)\in \cS_k^+\times \cS_k^-$, a trip arc connects the pick-up vertex $(s^+,s_k)$ to the drop-off vertex $(s^-,e_k)$. The subset of trip arcs associated with request $k$ is denoted by $\cA_{\texttt{trip}}^k$.

For each arc $a \in \cA_d$, let $\Delta t_a$ denote its duration, and let $\Delta b_a$ denote its net change in battery charge. Thus, $\Delta b_a$ is negative for trip arcs and nonnegative for waiting arcs. Let $\cA_v^-$ and $\cA_v^+$ denote the sets of arcs entering and leaving vertex $v$, respectively.

Let $f_a$ be the binary variable that equals 1 if arc $a \in \cA_d$ is used in the trajectory and 0 otherwise. Let $B$ denote the battery capacity of the homogeneous fleet, and let $SoC_v$ denote the state of charge at vertex $v \in \cV_d$, bounded in $[0,B]$. The pricing subproblem for representative day $d$ is then
\begingroup
\allowdisplaybreaks
\begin{subequations}
  \label{model:PP(MIP)}
\begin{align}
 \max \quad & \sum_{a \in \cA_d} \xi_a^d f_a - \mu_d \label{pp_obj} \\
\hbox{s.t.\ \ }
& \sum_{k \in \cK_d } \sum_{a \in \cA^k_{\texttt{trip}}} f_a \geq 1 \label{triptraj} \\
& \sum_{a \in \cA_{\bar{o}}} f_a = 1 \label{paths} \\
& \sum_{a \in \cA^-_{v}} f_a - \sum_{a \in \cA^+_{v}} f_a = 0 && \forall v \in \cV_d\setminus\{\bar{o},\bar{d}\} \label{flow_conservation} \\
& \sum_{a \in \cA_{\bar{d}}} f_a = 1 \label{pathe} \\
& SoC_v - SoC_u \leq \Delta b_a + B \cdot (1 - f_a) && \forall a=(u,v) \in \cA_d \label{soc_balance} \\
& SoC_{\bar{o}} = B \label{soc_initial} \\
& 0 \leq SoC_v \leq B && \forall v \in \cV_d \label{bounds_SoC} \\
& f_a \in \{0,1\} && \forall a \in \cA_d. \label{bound_arc}
\end{align}
\end{subequations}
\endgroup
Here, $\xi_a^d$ denotes the arc-level contribution induced by the dual multipliers in \eqref{eq:reduced_cost}. In particular, if $a \in \cA_{\texttt{trip}}^k$, then $\xi_a^d=-\alpha_k^d$; if $a$ is a waiting arc at station $s$, then $\xi_a^d$ equals the negative sum of the relevant charger-availability dual values over the discretized time periods occupied by that waiting arc; and all remaining arcs have zero contribution. Thus, the objective in \eqref{pp_obj} is exactly the reduced cost of a feasible trajectory on representative day $d$. 
Constraints \eqref{triptraj} ensure that each generated trajectory contains at least one trip. Constraints \eqref{paths}--\eqref{pathe} are flow-conservation constraints defining a source--sink path in the day-specific network. Constraints \eqref{soc_balance}--\eqref{bounds_SoC} enforce battery feasibility and state-of-charge bounds. Finally, constraints \eqref{bound_arc} specify the binary domain of the arc-selection variables.

Since $\cG_d$ is a directed acyclic graph, the pricing subproblem can be solved efficiently to optimality using a label-setting algorithm. Starting from the source node $\bar{o}$, the algorithm propagates labels forward through the network. For each vertex $v\in\cV_d$, it maintains a set of non-dominated labels, denoted by $\cL_v$, where each label represents a partial trajectory from $\bar{o}$ to $v$. A label stores the accumulated reduced contribution of the partial path together with the resource information needed for future extension decisions, most importantly the current state of charge. This pricing procedure is applied independently to each representative day, and its structure is unchanged across days.

Although the RMP is already feasible with no selected trajectories, we initialize the day-specific column pools $\{\cP_d'\}_{d\in\cD}$ with promising trajectories to accelerate convergence. To this end, we adapt the path-based heuristic of \citep{Brandstatter2020TS} separately to each representative day and use the resulting feasible trajectories as seed columns. Since joint compatibility with the shared station-opening, charger-installation, and fleet decisions is enforced by the master problem, these initial trajectories are generated day by day and need not be mutually compatible across different representative days.

%To start column generation, a feasible primal solution must be generated to serve as initial trajectory set $\cP'$. Such primal solutions can be constructed heuristically or by introducing artificial columns \citep{Desaulniers2006CGBook}. We adopt the path-based heuristic method introduced in \citep{Brandstatter2020TS} to generate such a solution. The algorithm starts with an empty path set, no stations are operational, no vehicles are purchased, and no trips are accepted. It simultaneously tracks the remaining budget $W'$ (initialized to the total budget $W$) and the number of vehicles $H_{s,t}$ entering station $s$ at time $t$. We identify a path with positive reduced cost. If such a path is found, the path set, remaining budget $W'$, and the number of vehicles $H_{s,t}$ entering stations are updated simultaneously.

%The arcs corresponding to all accepted trips are banned from $\cG_d$ by setting a large negative value. Meanwhile, for vertex $v=(s,t)$ where $H_{s,t}$ has reached the maximum capacity $Q_s$ at the station, all incoming arcs of these nodes are further banned by setting a large negative value, ensuring that paths selected subsequently can be integrated into the current solution without conflicts. The algorithm terminates when one of three scenarios is met: 1) the number of used vehicles reaches $H$; 2) the remaining budget is insufficient to support the purchase of a new vehicle; 3) no path that can generate additional revenue is detected.

\subsubsection{Exact branch-and-price algorithm}
Column generation solves only the LP relaxation of the selected single-objective model. Exact integer solutions are recovered by embedding column generation in a branch-and-bound scheme. Whenever the current RMP solution is fractional, branching is used to enforce integrality.

The branching rules fall into two groups. The first group directly affects the pricing problems because it restricts the structure of admissible trajectories. This includes branching on station-opening decisions ($y_s$), on trip pairs via Ryan--Foster branching, and on trip-acceptance variables ($x_k^d$). The second group acts only in the master problem. This includes branching on charger-installation variables ($z_s$) and, as a last resort, directly on trajectory variables ($\lambda_p^d$). We do not branch explicitly on $v_\texttt{h}$; once the trajectory decisions are integral, the cost-minimizing value of $v_\texttt{h}$ implied by \eqref{eq.vehicle_num} is integer automatically. Similarly, the fairness-related variables remain continuous and are determined by the master problem.

We apply the branching rules in a strict priority order. First, if any station-opening variable $y_s$ is fractional, we branch on a selected $y_s$ by enforcing either $y_s=0$ or $y_s=1$. The branch $y_s=0$ is enforced in the pricing problems by forbidding trajectories that visit station $s$ on any representative day. Second, if all $y_s$ variables are integral, we apply Ryan--Foster branching on a pair of trips from the same representative day that are fractionally assigned to the same vehicle, creating one branch in which the two trips must be served together and another in which they must be served apart. Third, if neither of the above applies, we branch on a fractional trip-acceptance variable $x_k^d$, yielding the two cases $x_k^d=0$ and $x_k^d=1$; the branch $x_k^d=0$ is enforced by forbidding all trip arcs associated with request $k$ in the pricing problem for day $d$.

If all variables above are integral but some charger-installation variables remain fractional, we branch on a selected $z_s$ using the standard split $z_s \le \lfloor \hat{z}_s \rfloor$ and $z_s \ge \lceil \hat{z}_s \rceil$. These constraints are handled in the master problem and affect pricing only indirectly through updated dual information. Finally, if all original variables are integral while some trajectory variables remain fractional, we branch directly on a column variable $\lambda_p^d$ by imposing $\lambda_p^d=0$ or $\lambda_p^d=1$. Although this rule is generally avoided, it guarantees completeness of the algorithm. Whenever a branching decision restricts the admissible trajectory structure, incompatible columns are removed from the current RMP and excluded from future pricing.

To explore the branch-and-bound tree, we use a hybrid node-selection strategy. Before any feasible integer solution is found, we adopt depth-first search to obtain incumbents quickly. Once a feasible solution is available, we switch to a best-bound strategy to accelerate convergence.

\subsection{Heuristic approach}
\label{sec.pareto_heuristic}
The exact balanced-box approach is computationally demanding because it may require many box-induced oracle calls, and each such oracle is itself solved by B\&P. For larger instances, this exact framework becomes prohibitively expensive. We therefore replace the exact B\&P oracle by a faster heuristic oracle based on diving.

Once the exact B\&P oracle is replaced by the diving heuristic, the balanced box structure from Section~\ref{sec.pareto_exact} is retained, but its interpretation changes. The anchor points are no longer exact anchors, and the boundary evaluations used to split boxes are no longer exact oracle values. As a consequence, the algorithm can no longer certify that a box is empty, nor can it guarantee that every nondominated point has been identified. The resulting method should therefore be regarded as a heuristic balanced box method. It preserves the same criterion-space logic as the exact framework, but the frontier it returns is an approximation of the Pareto frontier rather than an exact enumeration.

The diving heuristic operates after the column-generation process at the root node of the B\&P tree has converged. At this point, we have a fractional solution to the LP relaxation of the selected box-constrained oracle problem. Instead of initiating a full branch-and-bound search, the diving heuristic makes a sequence of deterministic fixing decisions to quickly construct a single integer-feasible solution. After each fixing, the RMP is re-solved, and the process continues along a single path of the conceptual branch-and-bound tree until an integer-feasible solution is obtained \citep{Sadykov2019IJoC}.

%In models (\textbf{M1})--(\textbf{M3}), 
In the box-constrained oracle problems, 
the discrete variables consist of shared strategic-design variables $(y_s,z_s,v_\texttt{h})$, day-specific trip-acceptance variables $x_k^d$, and day-specific trajectory-selection variables $\lambda_p^d$. By contrast, the fairness variables and their slack variables remain continuous and are updated automatically whenever the RMP is re-solved. In our preliminary computational experiments, directly applying a generic diving strategy to all discrete variables often failed to produce feasible solutions. We therefore implement a hierarchical fixing strategy that reflects the structure of the model.

Specifically, we adopt the fixing order $y_s \rightarrow \lambda_p^d \rightarrow z_s \rightarrow x_k^d$, while Ryan--Foster branching is ignored in the heuristic phase. The intuition is that station-opening decisions have the strongest structural impact on the admissible trajectory space, whereas trip-acceptance variables are the least fundamental once the strategic design and the trajectory pattern are largely fixed. We do not fix $v_\texttt{h}$ explicitly; once the trajectory decisions become integral, the smallest feasible value of $v_\texttt{h}$ is again determined automatically through \eqref{eq.vehicle_num}.

Within each variable class, the heuristic prioritizes values that are already close to integrality. For the binary variables $y_s$, $\lambda_p^d$, and $x_k^d$, we first search for fractional values close to 1 and tentatively fix them to 1. If no variable in the current class satisfies this criterion, we select the fractional variable with the largest value and again tentatively fix it to 1. If that fixing makes the RMP infeasible, we backtrack and fix the same variable to 0 instead. For station-opening decisions, we explicitly test both directions ($y_s=0$ and $y_s=1$) and keep the fixing that yields the better objective value after re-optimization.

For the general-integer charger-installation variables $z_s$, we fix a selected fractional variable to a nearby integer value. When both $\lfloor \hat z_s \rfloor$ and $\lceil \hat z_s \rceil$ are viable, we keep the one that yields the better objective value after re-optimizing the RMP. If one fixing direction leads to infeasibility, we backtrack and use the other. In this way, the diving heuristic remains consistent with the hierarchical structure of the exact B\&P algorithm while providing a much faster approximate oracle for large instances.

%%%%%%%%%%%%%%%%%%%%%%%%%%%%%%%%%%%%%%%%%%%%%%%%%%%%%%%%%%%%%%%%%%%%%%%%
%%%%%%%%%%%%%%%%%%%%%%%%%%%%%%%%%%%%%%%%%%%%%%%%%%%%%%%%%%%%%%%%%%%%%%%%
\section{Computational Experiments}\label{sec.experiments}
%%%%%%%%%%%%%%%%%%%%%%%%%%%%%%%%%%%%%%%%%%%%%%%%%%%%%%%%%%%%%%%%%%%%%%%%
%%%%%%%%%%%%%%%%%%%%%%%%%%%%%%%%%%%%%%%%%%%%%%%%%%%%%%%%%%%%%%%%%%%%%%%%
This section reports the computational results of the proposed methods. Section~\ref{sec.exp_setting} describes the instance families and parameter settings.
Section~\ref{sec.exp_artificial} evaluates the computational performance of the proposed solution framework.
Section~\ref{sec.exp_vienna_strategy} studies revenue--fairness trade-offs on the weekly Vienna  instances. Finally, Section~\ref{sec.exp_vienna_sensitivity} reports a sensitivity analysis on Vienna-inspired additionally generated instances and derives the corresponding managerial insights.

\subsection{Experimental settings and design}\label{sec.exp_setting}

We use three instance families, each serving a different purpose: synthetic grid-graph instances for scalability benchmarking against prior work, weekly Vienna instances for real-data fairness analysis, and Vienna-inspired  multi-day instances for controlled sensitivity analysis. All algorithms were implemented in C\#, and CPLEX 22.11 was used as the LP and MIP solver. The experiments were performed on a personal computer equipped with a 3.90 GHz Intel Core i9-12900K processor and 32 GB RAM.

\subsubsection{Benchmark settings}
\label{sec:benchmarksettings}
To rigorously evaluate the proposed models and algorithms, we use three instance families. 
The first family is used primarily for algorithmic scalability benchmarking against prior work, the second for real-data frontier and policy analysis, and the third for controlled sensitivity analysis under systematic variations in budget and demand composition.

The first family is generated following the grid-graph procedure in \citet{Brandstatter2020TS}. The operational area is modeled as a $50 \times 50$ grid graph, representing an urban area of approximately 96 km$^2$. Traversal times are set to 3 minutes for horizontal movements and 2 minutes for vertical movements. From the 2,500 vertices of this grid, different scenarios are created by varying the number of potential station locations. For each potential station, the strategic parameters $c_s$, $c'_s$, and $Q_s$ are generated as integers from $[9000,64000]$, $[22000,32000]$, and $[1,20]$, respectively, and the battery capacity is fixed at $B=100$. Trips are randomly generated from vertices located within a 5-minute travel radius of at least one station. Trips are divided into short-trip and long-trip groups. Battery consumption is drawn from $(5,25]$ for short trips and $(25,75]$ for long trips, trip duration is given by $d_k=\lceil \tau b_k\rceil$ with $\tau \sim U[3,4.5]$, and revenue is set to $p_k=30d_k$ cents. Unless otherwise stated, the vehicle cost is fixed at $c_\texttt{h}=20000$ and the total budget is fixed at $W=10{,}000{,}000$. We generate 12 benchmark instances for the single-representative-day setting by varying $|\cS|$, $|\cK|$, and $H$, where $H=|\cK|/5$ and $|\cK_{\texttt{low-energy}}|:|\cK_{\texttt{high-energy}}|=2:1$. To accommodate the multi-representative-day setting, we construct four additional instances that preserve the same total trip volume while distributing demand across multiple representative days. 
The instances are named according to the convention \texttt{D[\#days]S[\#stations]K[\#trips]H[\#vehicles]}, which provides a compact reference to their scale.

The second family consists of the Vienna benchmark introduced by  \cite{Brandstatter2020TS}, which provides a realistic urban network, candidate charging-station locations, and trip-request data for a station-based electric car-sharing system. We use four consecutive weeks of Vienna demand of the D7 operating area, for which the candidate-station sets contain 157--159 stations.
We construct four weekly multi-day instances, one for each observed week. Each weekly instance shares one strategic network and contains seven operating days. The resulting weekly trip totals are 520--538, with an average of 527.8.

To complement the real Vienna instances with a setting in which demand composition can be varied in a fully controlled manner, we construct a third family of Vienna-inspired multi-day instances.
We retain the D7 geographic setting and a fixed 50-station candidate network.
The trip pool is formed by combining de-duplicated observed D7 trips with additional constructed trips generated from observed origins, destinations, departure times, durations, and revenues, while battery consumption is recomputed from the corresponding network distance and travel time.
Each artificial instance contains five representative days with 400 trip requests per day, yielding 2000 trip requests in total. The fleet size is fixed at 80 vehicles so as to maintain the specification $H=|\cK|/5$ on a per-day basis, where $|\cK|=400$ refers to the number of trips on each representative day. To study the effect of market composition, we generate ten demand scenarios, consisting of five low-energy-dominant cases and five high-energy-dominant cases. These scenarios are later evaluated under multiple budget ratios. Their target low-energy shares are calibrated from the empirical D7 demand composition and range from 29.6\% to 70.4\%. 
%In total, we consider 40 Vienna-based  instances. 

\subsection{Algorithm performance evaluation}\label{sec.exp_artificial}
In this section, we evaluate the computational performance and scalability of the proposed solution framework.
First, in Section \ref{sec.exp_scalability}, we focus on the pure revenue-maximization variant of the problem and compare the performance of the state-of-the-art compact MIP model with the proposed B\&P algorithm on instances of different scales.
Specifically, the benchmark MIP model is derived from the FC formulation of \citep{Brandstatter2020TS}, which outperforms the alternative formulations considered there. 
Since our model is constructed for multiple representative days  whereas the FC formulation does not distinguish representative days, this benchmark comparison is restricted to the single-day setting.
Next, in Section \ref{sec.exp_diving}, we compare the B\&P algorithm with the diving heuristic and use the latter as a warm start to accelerate both the B\&P algorithm and the state-of-the-art MIP model.
Finally, in Section \ref{sec.exp_scalability_fair}, we test the performance of the proposed algorithms on the fairness-incorporated models in the multi-day setting.

\subsubsection{Exact methods scalability analysis}\label{sec.exp_scalability}
We first analyze the scalability of two exact methods for solving formulation (\textbf{M1}): the compact MIP formulation derived from \citep{Brandstatter2020TS} and our proposed B\&P algorithm. To ensure a controlled comparison, CPLEX is configured to use the traditional branch-and-cut search (\texttt{Cplex.MIPSearch.Traditional}) and a single computational thread. A uniform time limit of 7,200 seconds is imposed for all runs, and the time discretization level is initially set to 1 minute. Since the finest time granularity used in \citep{Brandstatter2020TS} is 15 minutes, we also report the corresponding comparison under 15-minute discretization in Table~\ref{tab:mip_vs_bap_15}.

The experiments are conducted on the eight small- and medium-scale instances. Table~\ref{tab:mip_vs_bap_results} reports the lower bound (LB), upper bound (UB), CPU time, optimality gap, and root-node LP objective value (Obj.Root).

\begin{table}[htbp]
  \centering
  \caption{Comparison of CPLEX (MIP) and the proposed Branch-and-Price (B\&P) algorithm.}
  \label{tab:mip_vs_bap_results}
  \resizebox{0.95\textwidth}{!}{
  \begin{tabular}{lrrrr rrrrr}
    \toprule
    \multirow{2}{*}{Instance Name} & \multicolumn{4}{c}{MIP~\citep{Brandstatter2020TS}} & \multicolumn{5}{c}{B\&P} \\
    \cmidrule(r){2-5} \cmidrule(l){6-10}
    & LB & UB & time (s) & gap (\%) & Obj.Root & LB & UB & time (s) & gap (\%) \\
    \midrule
    D1S10K50H10   & 1,463.40 & 1,463.40 & 0.01    & 0.00 & 1,463.40 & 1,463.40 & 1,463.40 & 0.02    & 0.00 \\
    D1S10K100H20  & 3,719.70 & 3,725.40 & 7,200.00 & 0.15 & 3,719.70 & 3,719.70 & 3,719.70 & 0.57    & 0.00 \\
    D1S10K200H40  & 7,794.30 & 8,010.27 & 7,200.00 & 2.77 & 7,939.80 & 7,939.80 & 7,939.80 & 21.40   & 0.00 \\
    D1S10K400H80  & 30.30    & 16,006.50 & 7,200.00 & -    & 15,798.30 & 14,685.60 & 15,798.30 & 7,200.00 & 7.58 \\
    D1S25K200H40  & 6,659.70 & 6,713.10 & 7,200.00 & 0.80 & 6,659.70 & 6,659.70 & 6,659.70 & 0.52    & 0.00 \\
    D1S25K400H80  & -        & 15,385.50 & 7,200.00 & -    & 14,610.60 & 14,610.60 & 14,610.60 & 31.69   & 0.00 \\
    D1S25K600H120 & -        & 25,542.30 & 7,200.00 & -    & 23,327.40 & 23,324.40 & 23,326.63 & 7,200.00 & 0.01 \\
    D1S25K800H160 & -        & 34,532.70 & 7,200.00 & -    & 31,309.20 & 31,192.80 & 31,309.20 & 7,200.00 & 0.37 \\
    \bottomrule
  \end{tabular}
  }
\end{table}

Table~\ref{tab:mip_vs_bap_results} shows that the compact MIP deteriorates quickly as instance size grows. While the MIP handles the smallest instance, it fails to prove optimality on all larger instances within the two-hour time limit and, for several medium-scale cases, fails to produce any feasible integer solution. In contrast, the B\&P algorithm proves optimality on five of the eight tested instances and terminates with very small residual gaps on two of the remaining three. Moreover, on some instances the final lower bound of the compact MIP remains below the optimal value proved by B\&P, indicating that the weakness of the compact formulation is not limited to primal incumbent quality but also extends to dual convergence. At the same time, the table shows that exact optimization remains challenging when the planning horizon is discretized at the 1-minute level, which substantially enlarges the underlying time-expanded networks and the resulting trajectory space.

Because the problem studied here is strategic rather than real-time operational, we also report results under a coarser 15-minute discretization, which remains meaningful at the planning level and matches the granularity used in \citep{Brandstatter2020TS}. In addition, Table~\ref{tab:mip_vs_bap_15} reports the lower and upper bounds obtained by CPLEX after 5 minutes, 1 hour, and 2 hours, so that incumbent progress as well as final performance can be compared directly.

\begin{table}[htbp]
  \centering
  \caption{Comparison of MIP and B\&P with 15-minute granularity.}
  \label{tab:mip_vs_bap_15}
  \resizebox{0.99\textwidth}{!}{
\begin{tabular}{lrrrrrrrrrrrrr}
\toprule
    \multirow{2}{*}{Instance Name} & \multicolumn{8}{c}{MIP~\citep{Brandstatter2020TS}}           & \multicolumn{5}{c}{B\&P}    \\
            \cmidrule(r){2-9} \cmidrule(l){10-14}
        & LB~(5 min) & LB~(1 hr) & LB~(2 hr) & UB~(5 min) & UB~(1 hr) & UB~(2 hr) & time~(s) & gap~(\%) & Obj.Root & LB       & UB       & time~(s) & gap~(\%) \\\hline
D1S10K50H10     & 1463.40   & -          & -          & 1463.40   & -          & -          & 0.02    & 0.00    & 1463.40  & 1463.40  & 1463.40  & 0.02    & 0.00    \\
D1S10K100H20     & 3691.80   & -          & -          & 3691.80   & -          & -          & 27.41   & 0.00    & 3691.80  & 3691.80  & 3691.80  & 0.17    & 0.00    \\
D1S10K200H40    & 498.90    & 7629.90    & 7701.90    & 7939.50   & 7915.20    & 7915.20    & 7200.00 & 2.77    & 7824.90  & 7824.90  & 7824.90  & 12.01   & 0.00    \\
D1S10K400H80    & 0.00      & 11352.00   & 13927.50   & 15845.70  & 15845.70   & 15845.70   & 7200.00 & 13.77   & 15614.70 & 14484.30 & 15614.70 & 7200.00 & 7.80    \\
D1S25K200H40    & 6447.90   & 6595.80    & 6595.80    & 6790.50   & 6681.30    & 6645.30    & 7200.00 & 0.75    & 6595.80  & 6595.80  & 6595.80  & 0.50    & 0.00    \\
D1S25K400H80    & 9159.90   & 13563.60   & 14034.60   & 15194.70  & 15194.70   & 15194.70   & 7200.00 & 8.27    & 14430.00 & 14430.00 & 14430.00 & 21.83   & 0.00    \\
D1S25K600H120   & 0.00      & 14231.40   & 18324.00   & 25542.30  & 23595.90   & 23595.90   & 7200.00 & 28.77   & 23146.80 & 21454.20 & 23146.80 & 7200.00 & 7.89    \\
D1S25K800H160   & -         & -          & -          & 34532.70  & 34532.70   & 34532.70   & 7200.00 & -       & 31204.80 & 28371.00 & 31204.80 & 7200.00 & 9.99   \\
\bottomrule
\end{tabular}
}
\end{table}

Table~\ref{tab:mip_vs_bap_15} confirms that coarser discretization reduces solution difficulty and enables the MIP to find feasible solutions on most instances.  
As expected, the objective values under 15-minute discretization differ from those under 1-minute discretization because the coarser grid changes the feasible trajectory set; these results are therefore intended as a benchmark at the literature granularity rather than as a substitute for the finer-resolution analysis provided in Table~\ref{tab:mip_vs_bap_results}. 
Even under this simplified setting, however, B\&P remains clearly superior in both runtime and solution quality. 
More specifically, the results show that, while CPLEX can obtain high-quality feasible solutions quickly on some of the easier instances, this behavior does not extend to the more difficult cases. On the harder instances, CPLEX sometimes improves the incumbent over time, but this improvement is often limited or eventually stalls, and sizable optimality gaps can remain even after two hours. More importantly, when both primal and dual progress are considered together, B\&P is stronger overall: it proves optimality on more cases, produces tighter bounds, and yields substantially smaller gaps on the challenging instances. In the largest case, CPLEX still fails to identify any feasible solution within the time limit. 

\subsubsection{Heuristic method evaluation}\label{sec.exp_diving}
To assess the value of diving as a feasible-solution oracle embedded in exact methods, we further examine its role in accelerating both B\&P and the direct MIP approach. 
To this end, we consider a hybrid two-phase strategy. In the first phase, the root node is solved by column generation, after which the diving heuristic is executed to obtain a high-quality integer solution. The integer solution found in this phase is used as the initial incumbent, thereby providing a strong upper bound at the start of the exact search. In the second phase, the full B\&P algorithm is launched with an RMP initialized by the enriched column pool collected from both the root-node solution and the diving procedure. 
For the MIP model, the integer solution is used as a MIP warm start. We denote these two variants as \texttt{Diving+B\&P} and \texttt{Diving+MIP}, respectively. 
All experiments in this subsection are conducted under the 15-minute discretization setting. Table~\ref{tab:diving_warm} reports the results for the standalone \texttt{Diving} heuristic, \texttt{Diving+MIP}, and \texttt{Diving+B\&P}. Since both \texttt{Diving} and \texttt{Diving+B\&P} attain the optimal solution on all eight instances, the table omits the UB, LB, and optimality-gap columns for these two methods and instead emphasizes their computational efficiency.

\begin{table}[!ht]
  \centering
  \caption{Comparison of diving heuristic and warm-started algorithms.}
  \label{tab:diving_warm}
  \resizebox{0.87\textwidth}{!}{
  \begin{tabular}{lrrrrrrr}
    \toprule
    && & \multicolumn{3}{c}{Diving+MIP} & Diving+B\&P & Diving \\
    \cmidrule(r){4-6} \cmidrule(l){7-7}\cmidrule(l){8-8}
    Instance Name & Obj.Root & LB & UB & time (s) & gap (\%) & time (s) & time (s) \\
    \midrule
    D1S10K50H10   & 1,463.40 & 1,463.40 & 1,463.40 & 0.12    & 0.00  & 0.04    & 0.04 \\
    D1S10K100H20   & 3,691.80 & 3,691.80 & 3,691.80 & 27.41   & 0.00  & 0.39    & 0.23 \\
    D1S10K200H40  & 7,824.90 & 7,824.90 & 7,939.50 & 7,200.00 & 1.46  & 3.29    & 1.14 \\
    D1S10K400H80  & 15,614.70 & 15,614.70 & 15,845.70 & 7,200.00 & 1.48  & 43.19   & 43.16 \\
    D1S25K200H40  & 6,595.80 & 6,595.80 & 6,698.40 & 7,200.00 & 1.56  & 0.23    & 0.23 \\
    D1S25K400H80  & 14,430.00 & 14,430.00 & 15,194.70 & 7,200.00 & 5.30  & 2.94    & 2.94 \\
    D1S25K600H120 & 23,146.80 & 23,146.80 & 23,595.90 & 7,200.00 & 1.94  & 66.07   & 66.07 \\
    D1S25K800H160 & 31,206.00 & 31,205.70 & 34,532.70 & 7,200.00 & 10.66 & 1,071.86 & 351.97 \\
    \bottomrule
  \end{tabular}
  }
\end{table}

Table~\ref{tab:diving_warm} shows that the diving heuristic is highly effective both as a standalone method and as a warm start for B\&P. The standalone heuristic produces optimal solutions very quickly on all instances, while \texttt{Diving+B\&P} proves optimality for all eight cases. By contrast, \texttt{Diving+MIP} still inherits the convergence difficulties of the compact MIP on larger instances. These results justify using the diving heuristic as the main practical solution tool once fairness constraints are introduced.

\subsubsection{Performance on fairness-incorporated formulations}\label{sec.exp_scalability_fair}
We test the performance of \texttt{Diving} and \texttt{Diving+B\&P} on formulations \textbf{M2} ($\epsilon^{\min}=0.7$) and \textbf{M3} ($\epsilon^{\text{gap}}=0.1$) under 15-minute discretization, with a 7200-second runtime limit for all 5-representative-day instances. The results are summarized in Table \ref{tab:diving_fair}.

\begin{table}[htbp]
  \centering
  \caption{Performance of diving heuristic and diving+B\&P on formulations (\textbf{M2})--(\textbf{M3}).}
  \label{tab:diving_fair}
  \resizebox{0.9\textwidth}{!}{
  \begin{tabular}{lrrrrrrrr}
    \toprule
    & \multicolumn{4}{c}{M2 ($\epsilon^{\min}=0.7$)} & \multicolumn{4}{c}{M3 ($\epsilon^{\text{gap}}=0.1$)} \\
    \cmidrule(r){2-5}\cmidrule(r){6-9}
    Instance Name & \multicolumn{2}{c}{Diving} & \multicolumn{2}{c}{Diving+B\&P} & \multicolumn{2}{c}{Diving} & \multicolumn{2}{c}{Diving+B\&P} \\
    \cmidrule(r){2-3}\cmidrule(r){4-5}\cmidrule(r){6-7}\cmidrule(r){8-9}
\multicolumn{1}{c}{}        & \multicolumn{1}{c}{time(s)} & \multicolumn{1}{c}{Gap(\%)} & \multicolumn{1}{c}{time(s)} & \multicolumn{1}{c}{Gap(\%)} & \multicolumn{1}{c}{time(s)} & \multicolumn{1}{c}{Gap(\%)} & \multicolumn{1}{c}{time(s)} & \multicolumn{1}{c}{Gap(\%)} \\\hline
D5S25K200H40  & 0.34 & 0.00 & 0.34 & 0.00 & 0.58 & 0.35 & 7200.00 & 0.17 \\
D5S25K400H80  & 0.39 & 0.00 & 0.39 & 0.00 & 0.25 & 0.03 & 7200.00 & 0.01 \\
D5S25K600H120 & 0.79 & 0.00 & 0.79 & 0.00 & 1.30 & 0.30 & 7200.00 & 0.30 \\
D5S25K800H160 & 3.95 & 0.72 & 7200.00 & 0.72 & 3.65 & 0.03 & 7200.00 & 0.03 \\
\bottomrule
  \end{tabular}
  }
\end{table}

Table~\ref{tab:diving_fair} indicates that fairness constraints increase the computational difficulty substantially, particularly under the service-rate disparity formulation. 
For the max-min model (\textbf{M2}), \texttt{Diving+B\&P}  proves optimality almost immediately on the three smaller instances and reaches the time limit only on the largest one. By contrast, for the service-rate disparity model (\textbf{M3}), \texttt{Diving+B\&P}  reaches the 7200-second time limit on all four instances and yields only marginal improvement over the initial heuristic incumbent. 
In contrast, the standalone diving heuristic consistently generates high-quality solutions with very small optimality gaps within 4 seconds across all instances, regardless of which fairness scheme is implemented. The heuristic shows exceptional robustness to both instance scale and constraint tightness. Therefore, for the substantive fairness analysis in the remainder of this section, the standalone diving heuristic is adopted as the practical solution method.

\subsection{Analysis of fairness strategies on Vienna instances}\label{sec.exp_vienna_strategy}
\begin{figure}[!ht]
    \centering
    \includegraphics[width=0.88\linewidth]{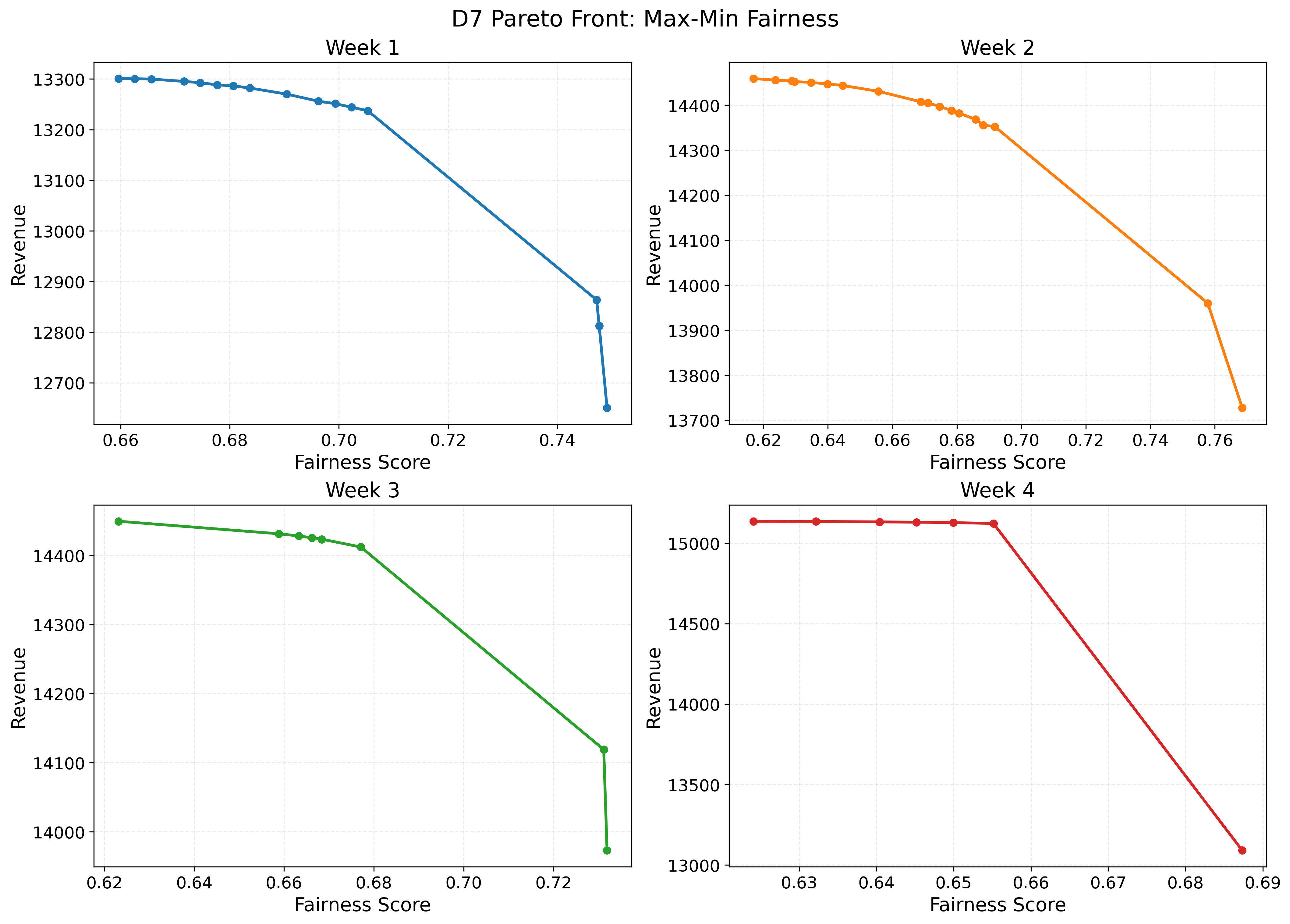}
    \caption{Pareto fronts of the max-min fairness model (\textbf{M2}) on the four D7 weekly instances.}
    \label{fig:vienna_maxmin_d7}
\end{figure}

\begin{figure}[!ht]
    \centering
    \includegraphics[width=0.88\linewidth]{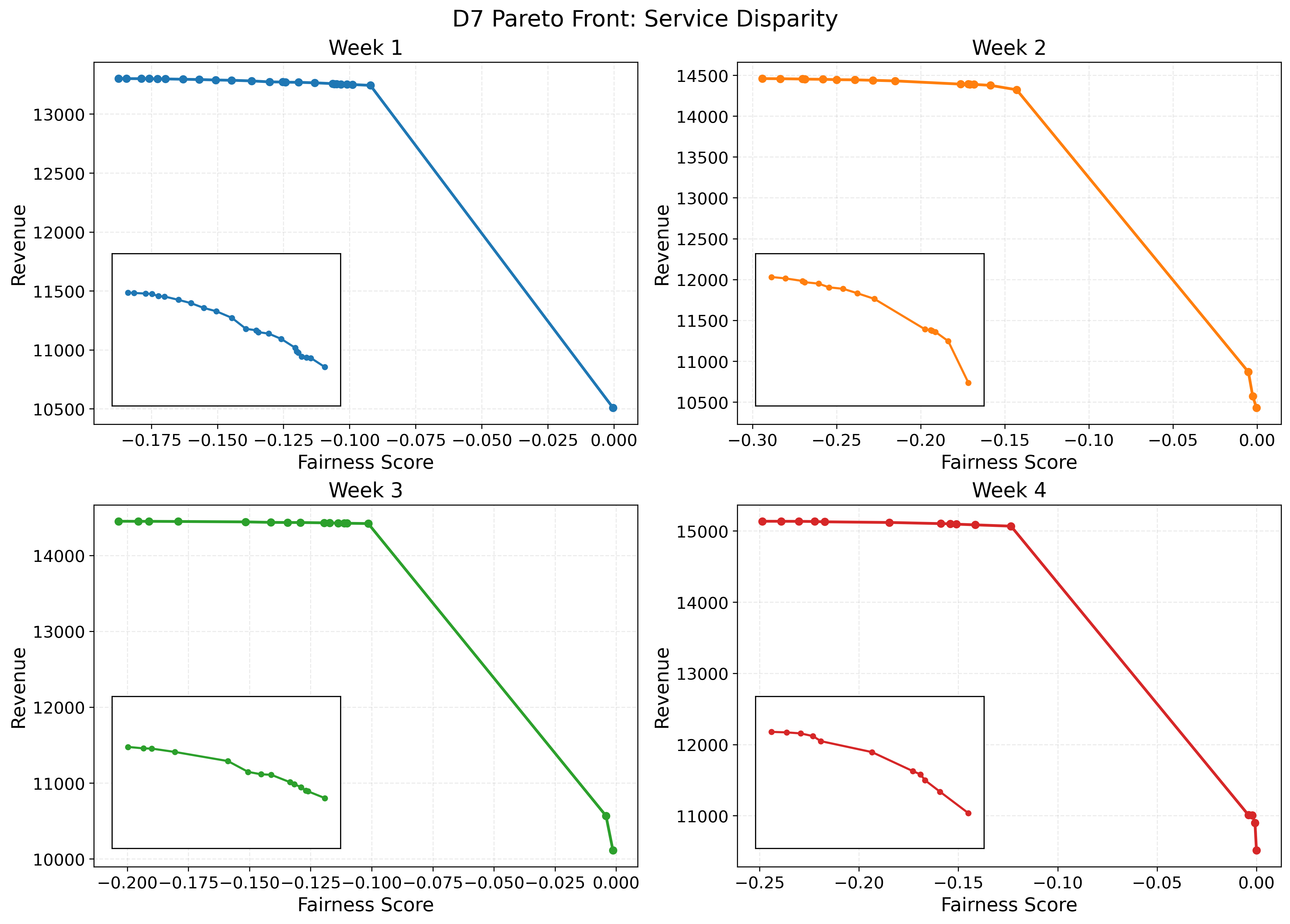}
    \caption{Pareto fronts of the service-rate disparity model (\textbf{M3}) on the four D7 weekly instances. }
    \label{fig:vienna_disparity_d7}
\end{figure}

We next examine the real-data revenue--fairness trade-offs on the weekly filtered Vienna instances. The analysis focuses on D7, the largest operating area among the three Vienna families and the one with the richest filtered demand among the Vienna instances in \citep{Brandstatter2020TS}. For each of the four D7 weekly instances, we generate Pareto-front approximations by the heuristic balanced box method described in Section~\ref{sec.pareto_heuristic}. For model (\textbf{M2}), the horizontal axis is the minimum average group service rate $\min_{g\in\cG}\bar{u}_g$. For model (\textbf{M3}), it is $-|\bar{u}_{\texttt{low}}-\bar{u}_{\texttt{high}}|$, so that both objectives remain in maximization form. Each Pareto point corresponds to one strategic design evaluated over seven representative operating days. The reported fronts therefore represent weekly multi-day operational performance rather than single-day snapshots.
Figures~\ref{fig:vienna_maxmin_d7} and \ref{fig:vienna_disparity_d7} show two distinct frontier geometries.  The horizontal axis of Figures~\ref{fig:vienna_maxmin_d7} reports the minimum average group service rate, whereas in Figure~\ref{fig:vienna_disparity_d7} it reports the negative absolute service-rate gap between the two groups; i.e., the respective fairness scores. The vertical axis in both figures reports expected revenue. For Figure~\ref{fig:vienna_disparity_d7}, the inset highlights the economically most relevant portion of the frontier before the steep revenue drop near exact parity.

Figure~\ref{fig:vienna_maxmin_d7} indicates that the max-min fronts are regular and gradually decreasing. Across the four weeks, the revenue-maximizing anchor occurs at minimum service rates between 0.617 and 0.660. In Weeks~1--3, increasing the minimum service rate to 0.732--0.768 reduces revenue by only about 476--731, indicating that substantial improvement in worst-group protection can be achieved within a limited efficiency loss. Week~4 is less favorable: the minimum service rate increases only from 0.624 to 0.687, while revenue decreases by about 2{,}046. Figure~\ref{fig:vienna_disparity_d7} exhibits a markedly different pattern. In every week, the initial portion of the front is nearly flat: from the revenue-maximizing anchor to the last point before the sharp break, revenue decreases by only about 30--137, while the service-rate gap narrows substantially. By contrast, the final movement toward exact or near-exact parity is associated with a pronounced revenue deterioration. Across the four weeks, the final step toward zero disparity induces an additional revenue loss of about 2{,}734--4{,}059. The inset emphasizes that most of the economically relevant variation occurs before this final cliff.

\begin{table}[!ht]
\centering
\small
\caption{Selected design metrics along the D7 Pareto fronts. Each entry reports the range across the four weekly instances. For (\textbf{M3}), the moderate point is the last point before the sharp revenue cliff.}
\label{tab:d7_pareto_design}
  \resizebox{0.95\textwidth}{!}{
\begin{tabular}{llccccc}
\toprule
Fairness formulation & Pareto point & Fairness metric & Revenue & Stations & Chargers & Trips \\
\midrule
\textbf{M2} & Anchor & $0.617$--$0.660$ & $13{,}301$--$15{,}137$ & $130$--$144$ & $164$--$174$ & $376$--$384$ \\
\textbf{M2} & High-protection endpoint & $0.687$--$0.768$ & $12{,}651$--$13{,}973$ & $140$--$146$ & $164$--$171$ & $363$--$394$ \\\hline
\textbf{M3} & Anchor & $0.187$--$0.294$  & $13{,}301$--$15{,}137$ & $130$--$144$ & $164$--$174$ & $376$--$384$ \\
\textbf{M3} & Moderate point & $0.092$--$0.143$  & $13{,}243$--$15{,}068$ & $133$--$146$ & $166$--$176$ & $378$--$385$ \\
\textbf{M3} & Near-parity endpoint & $0.000$--$0.001$  & $10{,}113$--$10{,}513$ & $124$--$137$ & $153$--$160$ & $273$--$313$ \\
\bottomrule
\end{tabular}
}
\end{table}

Turning to Table~\ref{tab:d7_pareto_design}, we summarize how these frontier patterns translate into strategic-design outcomes at selected Pareto points. For model (\textbf{M2}), the high-protection endpoint remains comparatively close to the anchor in terms of stations, chargers, and trips. This is consistent with the smooth trade-off visible in Figure~\ref{fig:vienna_maxmin_d7} and suggests that moderate max-min protection can often be achieved within a largely stable infrastructure configuration. For model (\textbf{M3}), the moderate point also remains close to the revenue-maximizing anchor, whereas the near-parity endpoint is associated with a much smaller served-trip volume and a substantially lower revenue level. The near-parity solutions also do not coincide with an expanded network; they often open fewer stations and install fewer chargers than the anchor solutions. This indicates that exact parity is obtained mainly by sacrificing profitable demand rather than by investing in additional capacity.

Overall, the D7 results show that moderate fairness targets define  the most relevant policy region for both formulations. Under max-min fairness, the trade-off remains gradual and operationally manageable over the observed range. Under service-rate disparity control, most fairness gains are realized before the steep revenue drop associated with near-parity solutions. These real-data observations motivate the calibrated sensitivity analysis in the next subsection around a moderate disparity target rather than an extreme parity requirement.

\subsection{Sensitivity analysis on Vienna-inspired instances}\label{sec.exp_vienna_sensitivity}
This subsection discusses the managerial implications of fairness constraints, budget levels, and demand composition on strategic system design. In contrast to the previous subsection, the focus here is on generalizable insights derived from controlled experiments. 
We report averages over the low-energy-dominant and high-energy-dominant Vienna-inspired scenario groups introduced in Section~\ref{sec:benchmarksettings}. 
Here, a budget ratio $\texttt{br} \in \{0.4,0.7,1.0\}$ means that the available budget is set to $\texttt{br} \cdot W$, where $W$ is the baseline budget used in the Vienna-inspired  instances. 

\begin{figure}[H]
    \centering
    \includegraphics[width=0.9\textwidth]{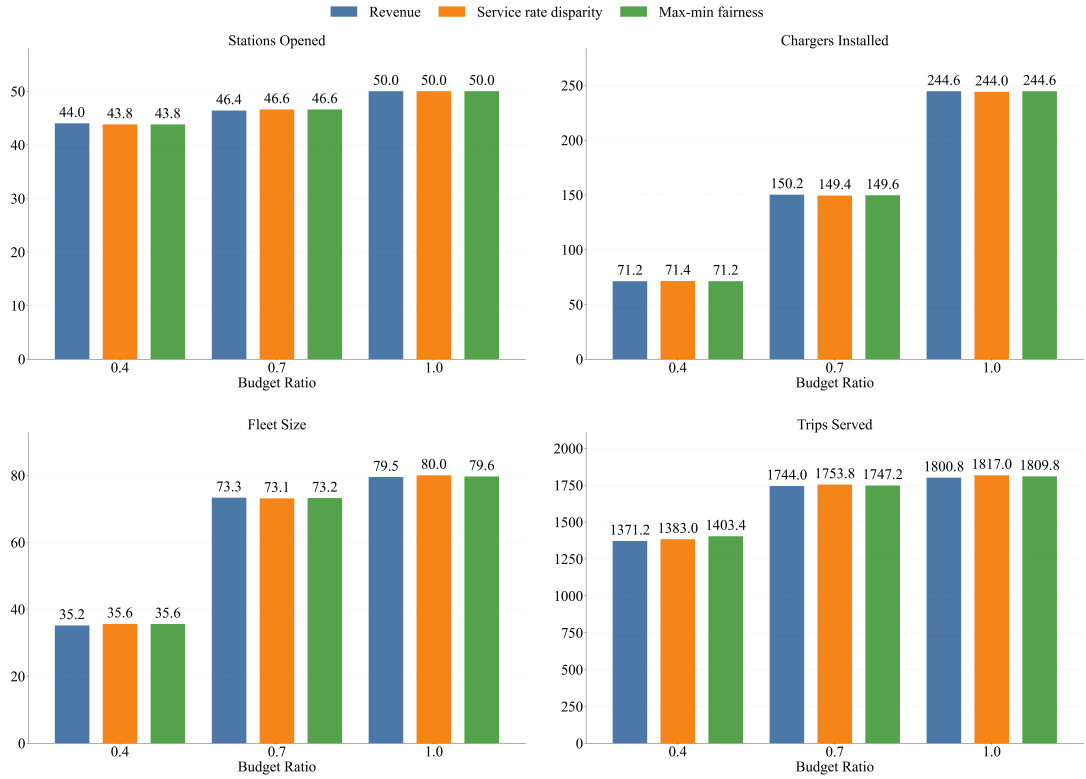}
    \caption{Impact of budget ratio on strategic decisions in low-energy-dominant scenario. }
    \label{fig:budget_impact_short}
\end{figure}

\begin{figure}[H]
    \centering
    \includegraphics[width=0.9\textwidth]{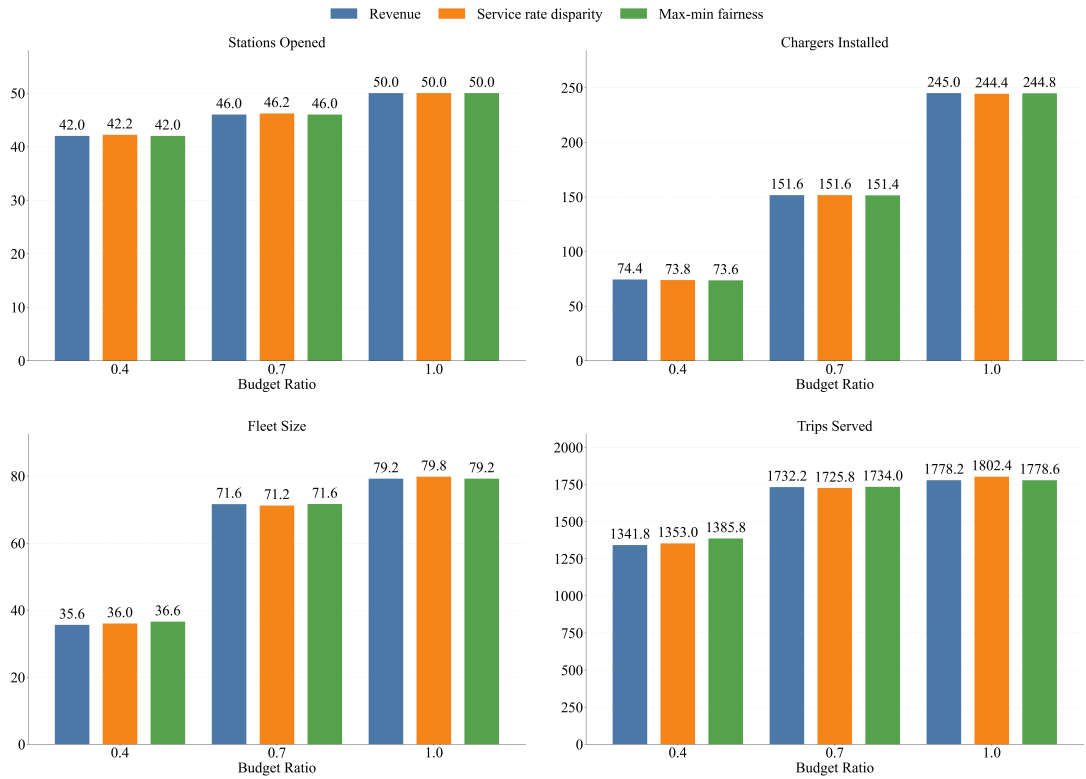}
    \caption{Impact of budget ratio on strategic decisions in the high-energy-dominant scenario.}
    \label{fig:budget_impact_long}
\end{figure}

\begin{figure}[H]
    \centering
    \includegraphics[width=0.9\textwidth]{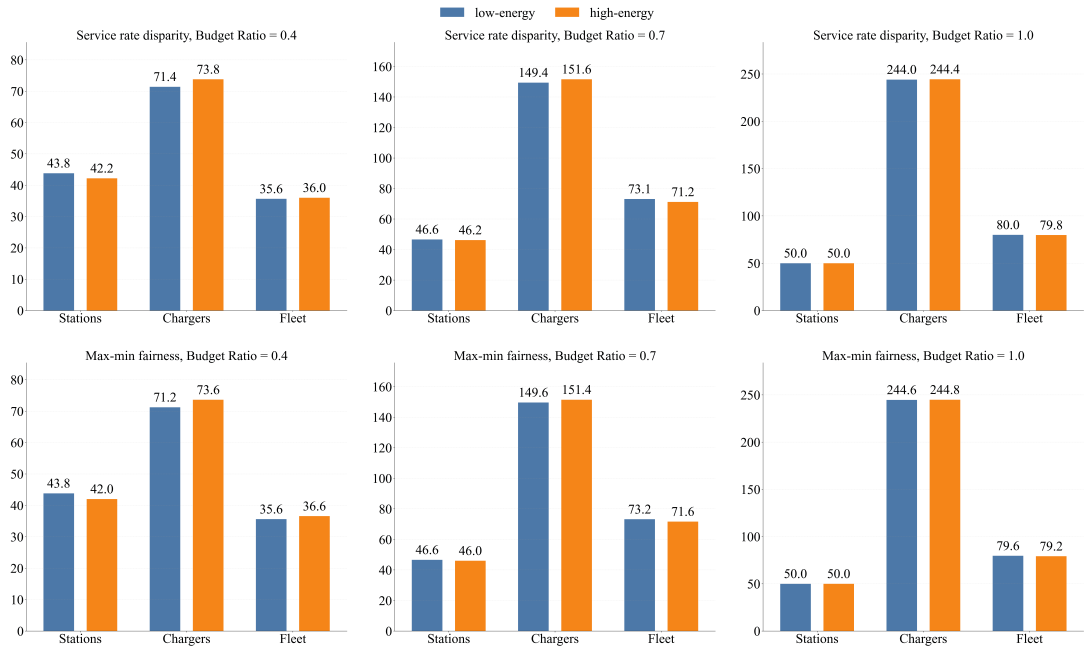}
    \caption{Impact of demand structure on strategic asset allocation.}
    \label{fig:demand_impact}
\end{figure}

\begin{figure}[H]
    \centering
    \includegraphics[width=0.9\textwidth]{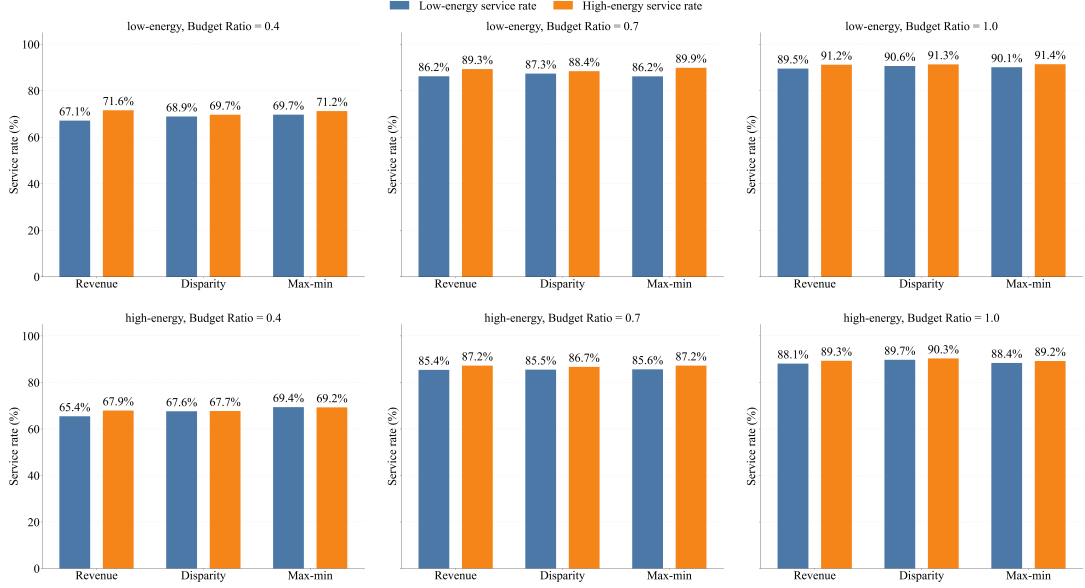}
    \caption{Impact of demand structure on service rate.}
    \label{fig:service_rate_impact}
\end{figure}

Figures~\ref{fig:budget_impact_short} and \ref{fig:budget_impact_long} present the average outcomes of three formulations---the revenue benchmark, service-rate disparity, and max-min fairness---across the three budget ratios in low-energy-dominant and high-energy-dominant markets, respectively, with respect to the numbers of stations opened, chargers installed, vehicles procured, and trips served. Figure~\ref{fig:demand_impact} further compares how demand composition affects strategic asset allocation under the two fairness models, and Figure~\ref{fig:service_rate_impact} reports the resulting group service rates. Taken together, these figures address two questions: how budget constraints reshape fairness-oriented strategic design, and how variations in demand composition alter the system’s asset portfolio.

Figures~\ref{fig:budget_impact_short} and \ref{fig:budget_impact_long} show that budget level does not merely scale the system uniformly, but largely determines whether fairness requirements can be supported by a sufficient operational-capacity base. As the budget ratio increases from 0.4 to 0.7, the number of opened stations rises only modestly, whereas the numbers of installed chargers, procured vehicles, and served trips all increase substantially. This indicates that the main improvement in this range comes from capacity expansion rather than network expansion. When the budget ratio further increases from 0.7 to 1.0, charger installation continues to grow markedly, but the increase in served trips becomes much smaller, suggesting that the system enters a regime of diminishing marginal service gains at higher budget levels. In other words, moving from low to medium budget primarily relaxes capacity bottlenecks, whereas moving from medium to high budget contributes more to capacity redundancy and operational flexibility than to proportional service expansion. Overall, these figures suggest that the main cost of pursuing fairness is concentrated in low-budget settings, where the system lacks the baseline capacity needed to support both efficiency and equity objectives.

Turning to demand composition, Figure~\ref{fig:demand_impact} shows that market structure primarily affects asset mix rather than network topology. Although both the low-energy-dominant and high-energy-dominant markets exhibit the same broad pattern of limited station growth and much faster growth in chargers and fleet size, the same fairness model induces different charger--fleet combinations across the two markets. This suggests that budget determines the overall scale boundary of the system, whereas demand composition determines how the available resources should be allocated within that boundary. Fairness-oriented strategic design therefore cannot be discussed independently of market structure: the same fairness objective may require different asset-allocation pathways under different demand compositions. Taken together, Figures~\ref{fig:budget_impact_short}--\ref{fig:demand_impact} show that budget, objective choice, and demand structure jointly shape the optimal strategic design of the system.

Compared with aggregate resource indicators, the service-rate results in Figure \ref{fig:service_rate_impact} provide a clearer view of how fairness constraints operate. Under the revenue benchmark (\textbf{M1}), one trip group systematically receives a lower service rate and therefore constitutes the disadvantaged group. The main role of the fairness-oriented formulations is not to substantially increase the total number of served trips, but to redistribute service quality across groups. In particular, the service-rate-disparity model (\textbf{M3}) more consistently compresses the inter-group service-rate gap, mainly by improving the lower-served group rather than by further strengthening the better-served group. Its fairness mechanism therefore acts more directly on the disparity objective and produces more stable effects across the tested demand scenarios. By contrast, the performance of the max-min model (\textbf{M2}) is more sensitive to the choice of the minimum-service target $\epsilon^{\min}$. If $\epsilon^{\min}$ is set too loosely, the resulting solution may not provide a meaningful improvement for the disadvantaged group; if it is set too aggressively, the model may require a substantial efficiency sacrifice before producing a visibly stronger convergence of the group service rates. This suggests that, although the max-min mechanism is appealing in principle, its practical behavior is more sensitive to calibration. Relative to (\textbf{M2}), formulation (\textbf{M3}) yields more stable and more interpretable fairness effects in the present experiments. At the same time, the service-rate differences across (\textbf{M1})--(\textbf{M3}) become much smaller as the budget increases, indicating that the tension between efficiency and fairness is concentrated mainly in low-budget settings.

Overall, these controlled experiments show that there is no universal fairness-oriented design blueprint. Budget levels determine whether the system can support the baseline capacity needed to pursue fairness objectives, whereas demand composition determines how the available capacity should be allocated between chargers and vehicles. Consequently, effective long-term planning for shared electric mobility requires not only the choice of a fairness policy, but also sufficient capital investment and a charger--fleet mix calibrated to the underlying market structure.

%%%%%%%%%%%%%%%%%%%%%%%%%%%%%%%%%%%%%%%%%%%%%%%%%%%%%%%%%%%%%%%%%%%%%%%%
%%%%%%%%%%%%%%%%%%%%%%%%%%%%%%%%%%%%%%%%%%%%%%%%%%%%%%%%%%%%%%%%%%%%%%%%
\section{Conclusion}
\label{sec.conclusion}
%%%%%%%%%%%%%%%%%%%%%%%%%%%%%%%%%%%%%%%%%%%%%%%%%%%%%%%%%%%%%%%%%%%%%%%%
%%%%%%%%%%%%%%%%%%%%%%%%%%%%%%%%%%%%%%%%%%%%%%%%%%%%%%%%%%%%%%%%%%%%%%%%
In this paper, we studied the strategic design of station-based one-way electric car-sharing systems under explicit fairness considerations. We developed an operations-aware planning framework that links long-term decisions on station opening, charger installation, and fleet sizing with day-level vehicle trajectories, charging feasibility, and realized service outcomes. To capture equity in a way that reflects actual system performance rather than static accessibility alone, we formulated and compared two fairness paradigms---service-rate disparity and max-min fairness---within a multi-day representative-demand setting. This led to a bi-objective trajectory-based model that jointly evaluates revenue and fairness across heterogeneous operating days. On the methodological side, we combined a tailored branch-and-price algorithm for the underlying oracle problems with an exact criterion-space search procedure to generate the Pareto frontier, and complemented this framework with a diving-based heuristic to obtain high-quality frontier approximations for larger instances. 

Our computational results demonstrate both the methodological value and the practical relevance of the proposed framework. Methodologically, the branch-and-price approach remains highly effective for the benchmark revenue-maximization setting studied in the literature, solving substantially larger and more finely time-discretized instances than alternative formulations. In the fairness-aware setting, the results show that the two fairness paradigms induce materially different trade-off structures and planning implications. On the weekly Vienna instances, max-min fairness yields comparatively smooth revenue--fairness trade-offs, whereas service-rate disparity exhibits a flat region of moderate improvement followed by a steep revenue drop near exact parity. On the Vienna-inspired  instances, we find that the main cost of pursuing fairness is concentrated in low-budget settings, where additional investment is used primarily to relax charger and fleet bottlenecks rather than to expand the station network. We also find that demand composition affects asset mix more than network topology, implying that fairness-oriented planning depends not only on the choice of fairness policy, but also on how charger and fleet resources are calibrated to the underlying market structure. Taken together, these findings show that fairness can be embedded directly into long-term electric car-sharing design in a computationally tractable way, while yielding actionable insights for planners and policymakers concerned with both economic viability and equitable service provision. 

Several directions remain for future research. One natural extension is to incorporate richer uncertainty by combining the present planning framework with stochastic demand models and behaviorally grounded user choice. Another is to enrich the operational layer with dynamic relocation and coordinated charging decisions, so that the interaction between fairness objectives and day-to-day fleet management can be studied more fully. It would also be valuable to consider adaptive or multi-stage investment planning, reflecting the gradual rollout of electric car-sharing systems in practice. Finally, the framework could be extended to compare alternative notions of fairness beyond the two studied here, thereby supporting a broader discussion of which equity concepts are most appropriate for different shared-mobility planning settings.

%\section*{Data availability}\vskip-6pt
%Data will be made available on request.

\section*{Acknowledgement}\vskip-6pt
The authors thank Georg Brandst\"{a}tter, Markus Leitner, and Ivana Ljubi\'{c} for sharing the Vienna dataset and for their assistance with the instance generation process.

%\vskip-18pt
%\newpage
\onehalfspacing
\bibliographystyle{elsarticle-harv}
\bibliography{main}

%\appendix
%\include{2025_8_Appendix}

\end{document}